\documentclass[11pt,reqno]{amsart}
\usepackage{mathrsfs}
\usepackage{amsmath}
\usepackage{amssymb, amsmath}
\usepackage{graphicx}
\numberwithin{equation}{section}

\def\F{\mathscr{F} }

\def\R{\mathbb{R}}

\def\T{\mathbb{T}}
\def\C{\mathbb{C}}

\newtheorem{thm}{Theorem}[section]
\newtheorem{lem}{Lemma}[section]

\newtheorem{prop}{Proposition}[section]

\newtheorem{cor}{Corollary}[section]

\newtheorem{defn}{Definition}[section]

\newtheorem{remark}{Remark}[section]

\makeatletter
\newcommand{\Extend}[5]{\ext@arrow0099{\arrowfill@#1#2#3}{#4}{#5}}
\makeatother

\begin{document}

\setcounter{page}{1}

\title[Global rough solution for DNLS]{Global well-posedness for Schr\"{o}dinger equation with derivative in $H^{\frac{1}{2}}(\R)$}
\author{Changxing Miao}
\address{Institute of Applied Physics and
Computational Mathematics,  P. O. Box 8009,\ Beijing,\ 100088, P. R.
China,\ } \email{miao\_changxing@iapcm.ac.cn}

\author{Yifei Wu}
\address{Department of Mathematics, South China University of Technology,
Guangzhou, Guangdong, 510640, P. R. China, }
\email{yerfmath@yahoo.cn}

\author{Guixiang Xu}
\address{Institute of
Applied Physics and Computational Mathematics,  P. O. Box 8009,\
Beijing,\ 100088, P. R. China, } \email{xu\_guixiang@iapcm.ac.cn}


\date{}


\maketitle

\begin{abstract}\noindent In this paper, we consider the Cauchy problem of the cubic nonlinear
Schr\"{o}dinger equation with derivative in $H^s(\R)$.
This equation was known to be the local well-posedness  for $s\geq
\frac12$ (Takaoka,1999), ill-posedness  for $s<\frac12$ (Biagioni
and Linares, 2001, etc.) and global well-posedness  for $s>\frac12$
(I-team, 2002). In this paper, we show that it is global
well-posedness  in $H^{\frac{1}{2}}(\R)$. The main approach is the
third generation I-method combined with some additional resonant
decomposition technique.  The resonant decomposition is applied to
control the singularity coming from the resonant interaction.
%
%
\end{abstract}

 \baselineskip=18pt

\section{Introduction}
In this paper, we consider the Cauchy problem of the Schr\"{o}dinger
equation with derivative:
 \begin{equation}\label{eq:DNLS}
   \left\{ \aligned
    i\partial_{t}u+\partial_{x}^{2}u=& i\lambda\partial_x(|u|^2u),\;\; x\in \R, t\in \R,
    \\
    u(0,x)  = & u_0(x)\in   H^s(\mathbb{R}),
   \endaligned
  \right.
 \end{equation}
where $\lambda\in \R$, $H^s(\R)$ denotes the usual inhomogeneous
Sobolev space of order $s$.  It arises from describing the
propagation of circularly polarized Alfv\'{e}n waves in the magnetized
plasma with a constant magnetic field (see \cite{MOMT-PHY, M-PHY,
SuSu-book}).

The local well-posedness for (\ref{eq:DNLS}) is well understood.
By the Fourier restriction norm in \cite{Bourgain-93-NLS,
Bourgain-93-KdV} and the gauge transformation  in \cite{Ha-93-DNLS,
HaOz-92-DNLS, HaOz-94-DNLS}, Takaoka obtained the local
well-posedness of (\ref{eq:DNLS}) in $H^s(\R)$ for $s\geq 1/2$ in
\cite{Ta-99-DNLS-LWP}. This result was shown by Biagioni and Linares
\cite{BiLi-01-Illposed-DNLS-BO}, Bourgain \cite{Bourgain-97-KdV} and
Takaoka \cite{Ta-01-DNLS-GWP} to be sharp in the sense that the flow
map fails to be uniformly $C^0$ for $s<1/2$.

The global well-posedness for (\ref{eq:DNLS}) was also widely
studied. In \cite{Oz-96-DNLS}, Ozawa made use of two gauge
transformations and the conservation of the Hamiltonian, and showed
that (\ref{eq:DNLS}) was globally well-posed in $H^1(\R)$ under the
condition (\ref{small condition}). In \cite{Ta-01-DNLS-GWP}, Takaoka
used Bourgain's ``Fourier truncation method"
(\cite{Bourgain-98-2D-NLS, Bourgain-book-NLS}) to obtain the
global well-posedness in $H^s(\R)$ for $s>\frac{32}{33}$, again
under (\ref{small condition}). In \cite{CKSTT-01-DNLS,
CKSTT-02-DNLS}, I-team (Colliander-Keel-Staffilani-Takaoka-Tao) made
use of the first, second generations of I-method
to obtain the global wellposedness in $H^s(\R)$, for $s>2/3$ and
$s>1/2$, respectively. For other results, please refer to
\cite{Gr-05-NLS-DNLS, Ha-93-DNLS, HaOz-92-DNLS, HaOz-94-DNLS,
He-06-DNLS-T, Oz-96-DNLS, TsFu-80-DNLS, TsFu-81-DNLS, Win-09-DNLS}.

In this paper, we will combine the third generation of the I-method
with the resonant decomposition
 to show the global well-posedness of \eqref{eq:DNLS} in  $H^{\frac{1}{2}}(\R)$. We think that the resonant
decomposition technique may also be used to study the global
well-posedness of \eqref{eq:DNLS} in $H^{\frac{1}{2}}(\T)$.

\begin{thm}
The Cauchy problem (\ref{eq:DNLS}) is globally well-posed in
$H^\frac{1}{2}(\R)$ under the assumption of
\begin{equation}\label{small condition}
\|u_0\|_{L^2}<\sqrt{\frac{2\pi}{|\lambda|}}.
\end{equation}
\end{thm}

The main approach, as described above, is the I-method. This method
is based on the correction analysis of some modified energies and an
iteration of local result. The first modified energy is defined as
$E(Iu)$, for some smoothed out operator $I$ (see (\ref{I})).
Moreover, one can effectively add a ``correction term'' to $E(Iu)$.
This gives the second modified energy $E_I^2(u)$, and allows us to
better capture the cancellations in the frequency space. However, a
further analogous procedure does not work again. Since in this
situation, a strong resonant interaction appears and this resonant
interaction will make the related multiplier to be singular. More
precisely, as shown in \cite{CKSTT-02-DNLS}, we define the second
modified energy by a 4-linear multiplier $M_4$, which will generate
a 6-linear multiplier $M_6$ in the increment of
the second modified energy. If we define the third modified energy
naturally by the 6-linear multiplier $\sigma_6$ as
$$
\sigma_6=-\frac{M_6}{\alpha_6},
$$
where
$\alpha_6=-i(\xi_1^2-\xi_2^2+\xi_3^2-\xi_4^2+\xi_5^2-\xi_6^2)$, then
$\alpha_6$ vanishes in some large sets but $M_6$ does not. So it is
not suitable to define the third modified energy in this way. Our
argument is to decompose the multiplier $M_6$ into two parts: one is
relatively small and another is non-resonant.
The analogous way of resonant decomposition
was previously used in \cite{Wu-Xu:2009:NLS-4,
Miao-Shao-Wu-Xu:2009:gKdV}. However, it is of great complexity here
and a dedicated multiplier analysis is needed in this situation. The
resonant decomposition technical was also appeared previously in \cite{BDGS,
Bourgain-04-NLS-T, CKSTT-08}. In particular, I-team \cite{CKSTT-08}
made use of  the second generation ``I-method'', a resonant
decomposition (in order to avoid the ``orthogonal resonant
interaction'') and an ``angularly refined bilinear Strichartz estimate"
to obtain the global well-posedness of mass-critical
nonlinear Schr\"{o}dinger equation in dimension two.

\begin{remark}
Without loss of generality, we may take $\lambda=1$ in (\ref{eq:DNLS}) in the following context. Indeed,
we may first assume that $\lambda> 0$, otherwise, we may consider $\bar{u}$ for instead.
Then we may rescale the solution by the
transformation
$$
u(x,t)\rightarrow\frac{1}{\sqrt{\lambda}}u(x,t).
$$
This deduces the general case to the case of $\lambda=1$.
\end{remark}

\begin{remark}
For the global well-posedness, it is natural to impose the condition
(\ref{small condition}). Indeed, the solution of (\ref{eq:DNLS})
(for $\lambda=1$) enjoys the mass and energy conservation laws
\begin{equation}\label{mass-law}
M(u(t)):=\int|u(t)|^2\,dx=M(u_0),
\end{equation}
and
\begin{equation}\label{energy-law}
H(u(t)):= \int\Big[| u_x(t)|^2+\frac{3}{2}\text{Im}
|u(t)|^2u(t)\overline{u_x(t)}+\frac{1}{2}|u(t)|^6\Big]\,dx=H(u_0).
\end{equation}
By a variant gauge transformation
$$
v(x,t):=e^{-\frac{i}{4}\int_{-\infty}^x |u(y,t)|^2\,dy}u(x,t),
$$
we have
\begin{equation*}
\aligned
\|v(t)\|_{L^2_x}&=\|u(t)\|_{L^2_x},\\
H(u(t))&=\|v_x(t)\|_{L^2_x}^2-\frac{1}{16}\|v(t)\|_{L^6_x}^6.
\endaligned
\end{equation*}
Thus, the condition (\ref{small condition}) guarantee the energy $H(u(t))$ to be positive
via the sharp Gagliardo-Nirenberg inequality
$$
\|f\|_{L^6}^6\leq \frac{4}{\pi^2}\|f\|_{L^2}^4\|f_x\|_{L^2}^2.
$$
\end{remark}

\begin{remark}
In \cite{CKSTT-01-DNLS}, I-team obtained the increment bound
$N^{-1+}$ of the first generation modified energy, which leads to
the global well-posedness in $H^s(\R)$ for $s>2/3$. In
\cite{CKSTT-02-DNLS}, the authors obtained the increment bound
$N^{-2+}$ of the second modified energy , which extend the exponent $s$ to $s>1/2$. In this paper, we will make use
of the resonant decomposition to show the  increment bound
$N^{-5/2+}$ of the third generation modified energy, which allows us
to extend the exponent $s$ to $s= 1/2$.
\end{remark}

The paper is organized as follows. In Section 2, we give some
notations and state some preliminary estimates that will be used
throughout this paper. In Section 3, we introduce the gauge
transformation and transform (\ref{eq:DNLS}) into another equation.
Then we present the conservation law and define the modified
energies. In Section 4, we establish the upper bound of the
multipliers generated in Section 3. In Section 5, we obtain an upper
bound on the increment of the third modified energy. In Section 6,
we prove a variant local well-posedness result. In Section 7, we
give a comparison between the first and third modified energy. In
Section 8, we prove the main result.

\section{Notations and Preliminary Estimates}

We use $A\lesssim B$, $B\gtrsim A$ or sometimes $A=O(B)$ to denote
the statement that $A\leq CB$ for some large constant $C$ which may
vary from line to line, and may depend on the data. When it is
necessary, we will write the constants by $C_1(\cdot),C_2(\cdot),\cdots$ to see
the dependency relationship. We use $A\sim B$ to mean $A\lesssim
B\lesssim A$. We use $A\ll B$, or sometimes $A=o(B)$ to denote the
statement $A\leq C^{-1}B$. The notation $a+$ denotes $a+\epsilon$
for any small $\epsilon$, and $a-$ for $a-\epsilon$.
$\langle\cdot\rangle=(1+|\cdot|^2)^{1/2}$,
$J_x^\alpha=(1-\partial^2_x)^{\alpha/2}$. We use
$\|f\|_{L^p_tL^q_x}$ to denote the mixed norm
$\Big(\displaystyle\int\|f(\cdot,t)\|_{L^q}^p\
dt\Big)^{\frac{1}{p}}$. Moreover, we denote $\F_x$ to be the Fourier
transformation corresponding to the variable $x$.

For $s,b\in \R $, we define the Bourgain space $X_{s,b}^{\pm}$ to be
the closure of the Schwartz class under the norm
\begin{equation}
\|u\|_{X_{s,b}^{\pm}}:= \left(\iint
\langle\xi\rangle^{2s}\langle\tau\pm\xi^2\rangle^{2b}|\hat{u}(\xi,\tau)|^2\,d\xi
d\tau\right)^{1/2},\label{X}
\end{equation}
and we write $X_{s,b}:= X_{s,b}^+$ in default. To study the
endpoint regularity,
we also need a slightly stronger space $Y_s^{\pm}$ (than
$X_{s,\frac{1}{2}}^{\pm}$),
\begin{equation}\label{working space} \|f\|_{Y_s^{\pm}} :=
\|f\|_{X_{s,\frac{1}{2}}^{\pm}}+\left\|\langle
\xi\rangle^s\hat{f}\right\|_{L^2_\xi L^1_\tau}.
\end{equation}
These spaces obey the embedding $Y_s^{\pm}\hookrightarrow
C(\R,H^s(\R))$. Again, we write $Y_s:= Y_s^+$. It
motivates the space $Z_s$ related to Duhamel term under the norm
\begin{equation}\label{Duhamel space}
\|f\|_{Z_s} := \|f\|_{X_{s,-\frac{1}{2}}}+\left\|\frac{\langle
\xi\rangle^s\hat{f}}{\langle\tau+\xi^2\rangle}\right\|_{L^2_\xi
L^1_\tau}.
\end{equation}

Let $s<1$ and $N\gg 1$ be fixed,  the Fourier multiplier operator
$I_{N,s}$ is defined as \begin{equation}
\widehat{I_{N,s}u}(\xi)=m_{N,s}(\xi)\hat{u}(\xi),\label{I}
\end{equation} where the multiplier $m_{N,s}(\xi)$ is a smooth,
monotone function satisfying $0<m_{N,s}(\xi)\leq 1$ and
 \begin{equation} m_{N,s}(\xi)=\biggl\{
\begin{array}{ll}
1,&|\xi|\leq N,\\
N^{1-s}|\xi|^{s-1},&|\xi|>2N.\label{m}
\end{array}
\end{equation} Sometimes we denote $I_{N,s}$ and $m_{N,s}$ as $I$ and $m$
respectively for short if there is no confusion.

It is obvious that the operator $I_{N,s}$ maps $H^s(\R)$ into
$H^1(\R)$  for any $s<1$. More precisely, there exists some positive
constant $C$ such that \begin{equation}
    C^{-1}\|u\|_{H^s}\leq \|I_{N,s}u\|_{H^1}
\leq
    CN^{1-s}\|u\|_{H^s}.\label{E-I}
\end{equation} Moreover, $I_{N,s}$ can be extended to a map (still denoted by
$I_{N,s}$) from $X_{s,b}$ to $X_{1,b}$, which satisfies that for any
$s<1,b \in \R$,
$$
    C^{-1}\|u\|_{X_{s,b}}\leq\|I_{N,s}u\|_{X_{1,b}}
\leq CN^{1-s}\|u\|_{X_{s,b}}.
$$

Now we recall some well-known estimates in the framework of Bourgain space (see \cite{CKSTT-02-DNLS},
for example). First, Strichartz's estimate gives us
\begin{equation}
              \|u\|_{L^6_{xt}}
   \lesssim
              \|u\|_{X_{0,\frac{1}{2}+}^{\pm}}.\label{XE1}
\end{equation}
This interpolates with the identity
$$
\|u\|_{L^2_{xt}}=\|u\|_{X_{0,0}},
$$
to give
\begin{equation}\label{XE2}
\|u\|_{L^q_{xt}}
   \lesssim
              \|u\|_{X_{0,\theta+}^{\pm}}, \mbox{ for } \theta\geq \frac{3}{2}(\frac{1}{2}-\frac{1}{q}).
\end{equation}
Moreover, we have
\begin{equation}
\|f\|_{L^\infty_x L^{\infty}_t} \lesssim
\|f\|_{Y_{\frac{1}{2}+}}.\label{XE4}
\end{equation}
Indeed, by Young's and Cauchy-Schwartz's inequalities, we have
$$
\|f\|_{L^\infty_{xt}}\leq \left\|\hat{f}\right\|_{L^1_\xi
L^1_\tau}\lesssim \left\|\langle
\xi\rangle^{\frac{1}{2}+}\hat{f}\right\|_{L^2_\xi L^1_\tau}.
$$
\begin{lem}
Let $f\in Y_s^{\pm}$ for any $s>0$, then we have
\begin{equation}
              \|f\|_{L^6_{xt}}
   \lesssim
              \|f\|_{Y_{s}^{\pm}}.\label{XE5}
\end{equation}
\end{lem}
\begin{proof}We only consider $Y_s$-norm. By the dyadic decomposition, we write
$f=\sum_{j=0}^\infty f_j$, for each dyadic component $f_j$ with the
frequency support $\langle\xi \rangle\sim 2^j$. Then, by (\ref{XE2})
and (\ref{XE4}), we have
 \begin{equation*}
 \aligned
\|f\|_{L^6_{xt}} & \leq \sum\limits_{j=0}^\infty
\|f_j\|_{L^6_{xt}}\leq \sum\limits_{j=0}^\infty
\|f_j\|_{L^q_{xt}}^\theta\|f_j\|_{L^\infty_{xt}}^{1-\theta}\\
&\leq \sum\limits_{j=0}^\infty
\|f_j\|_{X_{0,\frac{1}{2}}}^\theta\|f_j\|_{Y_\rho}^{1-\theta}
\lesssim \sum\limits_{j=0}^\infty
2^{\rho(1-\theta)j} \|f_j\|_{Y_0},
 \endaligned
 \end{equation*}
where $\rho>\frac{1}{2}$, and we choose $q=6-$ such that
$\theta=1-$. Choosing $q$ close enough to 6 such that
$s>\rho(1-\theta)$, then we have the conclusion by
Cauchy-Schwartz's inequality.
\end{proof}
Moreover, interpolating between (\ref{XE4}) and (\ref{XE5}), we have
\begin{equation}
              \|f\|_{L^q_{xt}}
   \lesssim
              \|f\|_{Y_{s_q}^{\pm}},\label{XE6}
\end{equation}
for any $q\in (6,+\infty)$ and $s_q>\frac{1}{2}(1-\frac{6}{q})$.

At last, we give some bilinear estimates. Define the Fourier
integral operators $I_{\pm}^s(f,g)$ by
\begin{equation}
\widehat{I_{\pm}^s(f,g)}(\xi,\tau)=\displaystyle\int_{\star}m_{\pm}(\xi_1,\xi_2)^s
\hat{f}(\xi_1,\tau_1)\hat{g}(\xi_2,\tau_2),\label{Is}
\end{equation}
where $\displaystyle\int_\star =\int_{\stackrel{\xi_1+\xi_2=\xi,}{
\tau_1+\tau_2=\tau}}\,d\xi_1 d\tau_1$, and
$$
m_{-}=|\xi_1-\xi_2|, \quad m_+=|\xi_1+\xi_2|.
$$
Then we have
\begin{lem}\label{Bi-linear}
For the Schwartz functions $f,g$, we have
\begin{align}
\left\|I_-^{\frac{1}{2}}(f,g)\right\|_{L^2_{xt}} \lesssim&\;
\|f\|_{X_{0,\frac{1}{2}+}^+}\,\|g\|_{X_{0,\frac{1}{2}+}^+},\label{2.14}\\
\left\|I_-^{\frac{1}{2}}(f,g)\right\|_{L^2_{xt}} \lesssim&\;
\|f\|_{X_{0,\frac{1}{2}+}^-}\,\|g\|_{X_{0,\frac{1}{2}+}^-},\label{2.15}\\
\left\|I_+^{\frac{1}{2}}(f,g)\right\|_{L^2_{xt}} \lesssim&\;
\|f\|_{X_{0,\frac{1}{2}+}^+}\,\|g\|_{X_{0,\frac{1}{2}+}^-}.\label{2.16}
\end{align}
\end{lem}
\begin{proof} See \cite{Wu-Xu:2009:NLS-4} for example.
\end{proof}

When $s=0$, by (\ref{XE2}) we have
\begin{equation}
\left\|I_{\pm}^{0}(f,g)\right\|_{L^2_{xt}}  \leq
\|f\|_{L^p_{xt}}\,\|g\|_{L^q_{xt}}\lesssim
\|f\|_{X_{0,b+}}\,\|g\|_{X_{0,b'+}}, \label{2.17}
\end{equation}
where
$$ \frac{1}{p}+\frac{1}{q}=\frac{1}{2},\quad
b=\frac{3}{2}\Big(\frac{1}{2}-\frac{1}{p}\Big),\quad
b'=\frac{3}{2}\Big(\frac{1}{2}-\frac{1}{q}\Big),$$
that is,
$b+b'=\frac{3}{4}$, and $b,b'\in
\big[\frac{1}{4},\frac{1}{2}\big]$.

Interpolating between the results in Lemma \ref{Bi-linear} and
(\ref{2.17}) twice, we have
\begin{cor}
Let $I^{s}_{\pm}$ be defined by (\ref{Is}), then for any $s\in
[0,\frac{1}{2}]$,
\begin{align*}
\left\|I_-^{s}(f,g)\right\|_{L^2_{xt}}  \lesssim& \;
\|f\|_{X_{0,b_1+}^+}\,\|g\|_{X_{0,b_2+}^+},\\
\left\|I_-^{s}(f,g)\right\|_{L^2_{xt}}  \lesssim& \;
\|f\|_{X_{0,b_1+}^-}\,\|g\|_{X_{0,b_2+}^-},\\
\left\|I_+^{s}(f,g)\right\|_{L^2_{xt}}  \lesssim&\;
\|f\|_{X_{0,b_1+}^+}\,\|g\|_{X_{0,b_2+}^-}.
\end{align*}
where $b_1=\frac{1}{2}(1-s'+s)$, $b_2=\frac{1}{4}(2s'+1)$ for any
$s'\in [s,\frac{1}{2}]$.
\end{cor}

In this paper, we just need the following crude estimates:
\begin{align}
\left\|I_-^{\frac{1}{2}-}(f,g)\right\|_{L^2_{xt}} \lesssim&\;
\|f\|_{X_{0,\frac{1}{2}-}^+}\,\|g\|_{X_{0,\frac{1}{2}-}^+},\label{2.18}\\
\left\|I_-^{\frac{1}{2}-}(f,g)\right\|_{L^2_{xt}} \lesssim&\;
\|f\|_{X_{0,\frac{1}{2}-}^-}\,\|g\|_{X_{0,\frac{1}{2}-}^-},\label{2.19}\\
\left\|I_+^{\frac{1}{2}-}(f,g)\right\|_{L^2_{xt}} \lesssim&\;
\|f\|_{X_{0,\frac{1}{2}-}^+}\,\|g\|_{X_{0,\frac{1}{2}-}^-}.\label{2.20}
\end{align}

Before the end of this section, we record the following forms of the mean
value theorem, which are taken from \cite{CKSTT-03-KDV}. To prepare
for it, we state a definition: Let $a$ and $b$ be two smooth
functions of real variables. We say that $a$ is controlled by $b$ if
$b$ is non-negative and satisfies $b(\xi)\sim b(\xi')$ for
$|\xi|\sim |\xi'|$ and
$$a(\xi) \lesssim b(\xi), \,\,a^\prime (\xi)
\lesssim \frac {b(\xi)}{|\xi|},\,\,a^{\prime\prime} \lesssim \frac
{b(\xi)}{|\xi|^2}.
$$
\begin{lem}\label{le:mvt}
If $a$ is controlled by $b$ and $|\eta|, |\lambda|\ll |\xi|$, then
we have
\begin{itemize}
\item (Mean value theorem)
\begin{equation}\label{MTV}
\left| a(\xi+\eta)-a(\xi)\right| \lesssim |\eta| \frac
{b(\xi)}{|\xi|}.
\end{equation}
\item (Double mean value theorem)
\begin{equation}\label{DMTV}
\left| a(\xi+\eta+\lambda)-a(\xi+\eta)-a(\xi+\lambda)+a(\xi)\right|
\lesssim |\eta||\lambda| \frac {b(\xi)}{|\xi|^2}.
\end{equation}
\end{itemize}
\end{lem}

\section{The Gause transformation, energy and the modified energies}
%

\subsection{Gauge transformation and conservation laws} First, we summarize some results presented in \cite{CKSTT-01-DNLS, CKSTT-02-DNLS}. We start by recalling the gauge transformation used in
\cite{Oz-96-DNLS} to improve the derivative nonlinearity presented in (\ref{eq:DNLS}).
\begin{defn}
We define the non-linear map $\mathscr{G}: L^2(\R)\rightarrow
L^2(\R)$ by
$$
\mathscr{G}f(x):= e^{-i\int_{-\infty}^x |f(y)|^2\,dy}f(x).
$$
The inverse transformation $\mathscr{G}^{-1}f$ is then given by
$$
\mathscr{G}^{-1}f(x):= e^{i\int_{-\infty}^x |f(y)|^2\,dy}f(x).
$$
\end{defn}

Set $w_0:=\mathscr{G}u_0$ and $w(t):=\mathscr{G}u(t)$ for all time
$t$. Then (\ref{eq:DNLS}) is transformed to
 \begin{equation}\label{eq:DNLS2}
   \Biggl\{ \aligned
    &i\partial_{t}w+\partial_{x}^{2}w=-iw^2\partial_x\bar{w}-\frac{1}{2}|w|^4w,
    \;\; w:\R\times [0,T]\mapsto \C,\\
    &w(0,x)  =w_0(x),\quad x\in \R,t\in \R.
   \endaligned
 \end{equation}
In addition, the smallness condition (\ref{small condition}) becomes
\begin{equation}\label{small condition2}
\|w_0\|_{L^2}<\sqrt{2\pi}.
\end{equation}

Note that the transform $\mathscr{G}$ is a bicontinuous map from
$H^s(\R)$ to itself for any $s\in [0,1]$, thus the global
well-posedness of (\ref{eq:DNLS}) is equivalent to that of
(\ref{eq:DNLS2}). Therefore,  from now on, we focus our attention to
(\ref{eq:DNLS2}) under the assumption (\ref{small condition2}).

\begin{remark}
For the equation without the derivative term in (\ref{eq:DNLS2}) (it is just the focusing, mass-critical Schr\"{o}dinger equation):
 \begin{equation*}
   \Biggl\{ \aligned
    &i\partial_{t}w+\partial_{x}^{2}w=-|w|^4w,
    \\
    &w(0,x)  =w_0(x),\quad x\in \R,t\in \R,
   \endaligned
 \end{equation*}
 it is global well-posedness below $H^{\frac{1}{2}}(\R)$ with the mass less than that of the ground state. Indeed, in \cite{Wu-Xu:2009:NLS-4},
 the authors
 proved that it is global well-posedness in $H^s(\R)$ for $s>\frac{2}{5}$. So the difficulty of the equation (\ref{eq:DNLS2})
comes mainly from the derivative term.
\end{remark}
\begin{defn}
For any $f\in H^1(\R)$, we define the mass by
$$
M(f)=\int|f|^2\,dx,
$$
and the energy $E(f)$ by
$$
E(f):= \int|\partial_xf|^2\,dx-\frac{1}{2}\text{Im}\int
|f|^2f\partial_x\bar{f}\,dx.
$$
\end{defn}
By the gauge transformation and the sharp Gagliardo-Nirenberg
inequlity, we have (see \cite{CKSTT-01-DNLS} for details)
\begin{equation}\label{GN}
\|\partial_xf\|_{L^2}\leq C(\|f\|_{L^2})E(f)^{\frac{1}{2}},
\end{equation}
for any $f\in H^1(\R)$ such that $\|f\|_{L^2}<\sqrt{2\pi}$.

Moreover, the solution of (\ref{eq:DNLS2}) obeys the mass and energy conservation
laws (see cf. \cite{Oz-96-DNLS}):
\begin{equation}\label{mass-energy law}
M(w(t))=M(w_0),\quad E(w(t))=E(w_0).
\end{equation}

\subsection{Definition of $n$-linear functional}
Let $w$ be the solution of (\ref{eq:DNLS2}) throughout the
following contents. For an even integer $n$ and a given function
$M_n(\xi_1,\cdots,\xi_n)$ defined on the hyperplane
\begin{equation}\label{Gamma_n}
\Gamma_n=\left\{(\xi_1,\cdots, \xi_n):\xi_1+\cdots+\xi_n=0\right\},
\end{equation}
we define the quantity
\begin{align}\label{Lambda}
  \Lambda_n(M_n;w(t))
 :=\int_{\Gamma_n}M_n(\xi_1,\cdots,\xi_n)
   \F_xw(\xi_1,t)\overline{\F_xw}(-\xi_2,t)\\
  \qquad\qquad \cdots  \F_xw(\xi_{n-1},t)\overline{\F_xw}(-\xi_n,t)\,d\xi_1\cdots d\xi_{n-1}.\nonumber
\end{align}
Then by (\ref{eq:DNLS2}) and a directly computation, we have
\begin{align}\label{dLMD}
\frac{d}{dt}\Lambda_n(M_n;w(t))
=&\Lambda_n(M_n\alpha_n;w(t))\\
 &
-i\Lambda_{n+2}\Big(\sum\limits_{j=1}^n X_j^2(M_n)\xi_{j+1};w(t)\Big)\nonumber\\
 & +\frac{i}{2}\Lambda_{n+4}\Big(\sum\limits_{j=1}^n (-1)^{j+1}
X_j^4(M_n);w(t)\Big) ,\nonumber
\end{align}
where
$$
\alpha_n=i\sum\limits_{j=1}^n (-1)^j\xi_j^2,
$$
and
$$
X_j^l(M_n)=M_n(\xi_1,\cdots,\xi_{j-1},\xi_j+\cdots+\xi_{j+l},\xi_{j+l+1},\cdots,
\xi_{n+l}).
$$
Observe that if the multiplier $M_n$ is invariant under the
permutations of the even $ \xi_j$ indices, or of the odd $ \xi_j$
indices, then so is the functional $\Lambda_n(M_n;w(t))$.

\noindent \textbf{Notations:} In the following, we shall often write $\xi_{ij}$ for $\xi_i+\xi_j$,
$\xi_{ijk}$ for $\xi_i+\xi_j+\xi_k$, etc.. Also we write $m(\xi_i)=m_i$ and  $m(\xi_i+\xi_j)=m_{ij}$, etc..

\subsection{Modified Energies}
Define the first modified energy as
\begin{align}\label{E1}
E^1_I(w(t)):=& E(Iw(t))\\
=&-\Lambda_2\big(\xi_1\xi_2m_1m_2;w(t)\big)+\frac{1}{4}\Lambda_4\big(\xi_{13}m_1m_2m_3m_4;w(t)\big),\notag
\end{align}
where we have used the Plancherel  identity and (\ref{Lambda}).

We define the second
modified energy as
\begin{equation}\label{E2}
E^2_I(w(t)):=-\Lambda_2\big(\xi_1\xi_2m_1m_2;w(t)\big)+\frac{1}{2}\Lambda_4(M_4(\xi_1,\xi_2,\xi_3,\xi_4);w(t)),
\end{equation}
where
\begin{equation}\label{M4}
    M_4(\xi_1,\xi_2,\xi_3,\xi_4)=
    -\frac{m_1^2\xi_1^2\xi_3+m_2^2\xi_2^2\xi_4+m_3^2\xi_3^2\xi_1+m_4^2\xi_4^2\xi_2}
    {\xi_1^2-\xi_2^2+\xi_3^2-\xi_4^2}.
\end{equation}
Then by (\ref{dLMD}) (or see \cite{CKSTT-02-DNLS} for more
details), we have
\begin{equation}\label{dE2}
\frac{d}{dt}E^2_I(w(t))= \Lambda_6(M_6;w(t))
+\Lambda_{8}(M_{8};w(t)),
\end{equation}
where
\begin{align}\label{M6}
M_6(\xi_1,&\cdots,\xi_6) \\
  :=& \beta_6(\xi_1,\cdots,\xi_6)\notag\\
&-\frac{i}{72}
\sum\limits_{{\{a,c,e\}=\{1,3,5\}}\atop{\{b,d,f\}=\{2,4,6\}}}\;
\Big(M_4(\xi_{abc},\xi_d,\xi_e,\xi_f)\xi_b+M_4(\xi_a,\xi_{bcd},\xi_e,\xi_f)\xi_c
\nonumber\\
 &  +
M_4(\xi_a,\xi_b,\xi_{cde},\xi_f)\xi_d+M_4(\xi_a,\xi_b,\xi_c,\xi_{def})\xi_e\Big),\nonumber
\\
\label{M8}M_{8}(\xi_1, & \cdots,\xi_{8})  \\
:=& C_8\sum\limits_{{\{a,c,e,g\}=\{1,3,5,7\}}\atop{
\{b,d,f,h\}=\{2,4,6,8\}}}\!\!
\Big(M_4(\xi_{abcde},\xi_f,\xi_g,\xi_h)+
M_4(\xi_a,\xi_b,\xi_{cdefg},\xi_h) \notag \\
 &
\qquad \qquad \qquad \quad
-M_4(\xi_a,\xi_{bcdef},\xi_g,\xi_h)-M_4(\xi_a,\xi_b,\xi_c,\xi_{defgh})\Big)\nonumber
\end{align}
for some constant $C_8$ and
\begin{equation}\label{beta-6}
\beta_6(\xi_1,\cdots,\xi_6):=-\frac{i}{6} \sum\limits_{j=1}^6 (-1)^j m_j^2\xi_j^2.
\end{equation}
Note that $M_4,M_6,M_8$ are invariant under the
permutations of the even $ \xi_j$ indices, or of the odd $ \xi_j$
indices.

In order to consider the endpoint case, we also need to define the
third modified energy. Before constructing it, we shall do some
preparations. We adopt the notations that
$$
|\xi_1^{*}|\geq |\xi_2^{*}|\geq\cdots\geq |\xi_{6}^{*}|\geq\cdots\geq |\xi_{n}^{*}|.
$$
Moreover,
by the symmetry of $M_{6},M_8$ (and other multipliers defined later), we may restrict in $\Gamma_{n}$
(defined in (\ref{Gamma_n})) that
$$
|\xi_1|\geq |\xi_3|\geq\cdots\geq |\xi_{n-1}|,\quad |\xi_2|\geq |
\xi_4|\geq\cdots\geq |\xi_{n}|.
$$
Now we denote the sets
\begin{equation*}\aligned
\Upsilon=&\; \{(\xi_1,\cdots, \xi_{6})\in \Gamma_{6}: |\xi_1^*|\sim|
\xi^*_2|\gtrsim N\}, \\
\Omega_1=&\; \{(\xi_1,\cdots, \xi_{6})\in \Upsilon: |\xi_1|\sim |\xi_3|\gg |\xi_3^*| \mbox{ or } |\xi_2|\sim |\xi_4|\gg |\xi_3^*| \},\\
\Omega_2=&\; \{(\xi_1,\cdots, \xi_{6})\in \Upsilon: |\xi_1|\sim |\xi_2|\gtrsim
N\gg |\xi_3^*|, |\xi_1|^{\frac{1}{2}}|\xi_1+\xi_2|\gg |\xi_3^*|^{\frac{3}{2}} \},\\
\Omega_3=&\; \{(\xi_1,\cdots, \xi_{6})\in \Upsilon:
|\xi_3^{*}|\gg |\xi_4^{*}|\},
\endaligned\end{equation*}
and let
$$
\Omega=\Omega_1\cup\Omega_2\cup\Omega_3.
$$
\begin{remark}
In generally, $|M_6|$ is not controlled by $|\alpha_6|$,
this is the main difficulty lied in our problem. However, we exactly have
(see Lemma 4.9 for the proof)
$$
|M_6|\lesssim |\alpha_6|, \mbox{ for any } (\xi_1,\cdots, \xi_{6})\in \Omega.
$$
For this reason, $\Omega$ is referred to the non-resonant set.
\end{remark}

 Rewrite (\ref{dE2}) by
\begin{equation}\label{dE2'}
\frac{d}{dt}E^2_I(w(t))=\Lambda_{6}(M_6\cdot \chi_{\Gamma_6\backslash
\Omega};w(t)) +\Lambda_{6}(M_6\cdot
\chi_{\Omega};w(t))+\Lambda_{8}(M_{8};w(t)).
\end{equation}
Now we are ready to define the third modified energy $E^3_I(w(t))$. Let
\begin{equation}
E^3_I(w(t))=\Lambda_6(\sigma_6;w(t))+E^2_I(w(t)),\;\;
\sigma_6=-\frac{M_6}{\alpha_6}\cdot \chi_{\Omega}.\label{E3}
\end{equation}
Then by (\ref{dLMD}) and (\ref{dE2'}),  one has
\begin{equation}
\frac{d}{dt}E^3_I(w(t))=\Lambda_{6}(M_6\cdot
\chi_{\Gamma_6\backslash\Omega};w(t))+\Lambda_{8}(M_{8}+
\widetilde{M}_{8};w(t)) +\Lambda_{10}(M_{10};w(t)),\label{dE3}
\end{equation}
where $M_6,M_8$ defined in (\ref{M6}), (\ref{M8}) respectively, and

\begin{align}
\widetilde{M}_{8}=
&-i\sum\limits_{j=1}^6 X_j^2(\sigma_6)\xi_{j+1},\label{M8'}\\
M_{10}= &\frac{i}{2} \sum\limits_{j=1}^{6} (-1)^{j+1}
X_j^4(\sigma_6).\label{M10}
\end{align}

\begin{remark} By the dyadic decomposition, we restrict that
$$
|\xi_j^*|\sim N_j^*, \mbox{ for any } j=1, 2, \cdots .
$$
Now we give some explanations about the construction of $\Omega_j$.
We keep in mind the denominator of $\sigma_6$,
$$
\alpha_6=-i(\xi_1^2-\xi_2^2+\xi_3^2-\xi_4^2+\xi_5^2-\xi_6^2).
$$
On one hand, for the non-resonant region, we expect $|\alpha_6|$ has a
large lower bound in $\Omega$. On the other hand, we expect that the multipliers
$M_6$ has a small upper bound on the resonant region
$\Gamma_6\backslash\Omega$.
\begin{enumerate}
\item[\rm (a)] By the definition of $\Omega_1$, we have
$$
|\alpha_6|\sim {N_1^*}^2, \,\,\mbox{ for } (\xi_1,\cdots,\xi_6)\in \Omega_1.
$$
On the other hand, in $\Gamma_6\backslash\Omega_1$, the following case is ruled
out:
$$
\xi_1^*=\xi_1, \,\,\xi_2^*=\xi_3;\quad \mbox{ or }\quad \xi_1^*=\xi_2, \,\,\xi_2^*=\xi_4.
$$
Therefore, to estimate $M_6\cdot\chi_{\Gamma_6\backslash \Omega}$, we
only need to consider
$$
\xi_1^*=\xi_1, \,\,\xi_2^*=\xi_2;\quad \mbox{ or }\quad \xi_1^*=\xi_2, \,\,\xi_2^*=\xi_1.
$$
This is carried out in Proposition 4.1 below.

\item[\rm (b)] Now assume that we are in the situation: $|\xi_1|\sim
|\xi_2|\gtrsim N\gg |\xi_3^*|$. We find that $\alpha_6$ will not
vanish if
$$
|\xi_1+\xi_2|\gg {N_3^*}^2/N_1^*,
$$
since in this case $|\alpha_6|\sim |\xi_1||\xi_1+\xi_2|.$ It is
common to choose a lower bound of $|\xi_1+\xi_2|$ between
${N_3^*}^2/N_1^*$ and $N_3^*$, and the choice of the bound will
affect the bound of $M_6$ and $\widetilde{M}_8$. Generally (but not
absolutely), a small lower bound of $|\xi_1+\xi_2|$ gives a small
upper bound of $M_6$, but it maybe lead to a large upper bound of
$\widetilde{M}_8$. So, it appears important to make a suitable
choice.

As shown in the definition of $\Omega_2$, we choose a middle bound of
$$
|\xi_1+\xi_2|\gg |\xi_3^*|^{\frac{3}{2}}/|\xi_1|^{\frac{1}{2}}.
$$
This leads to the upper bound of $M_6\cdot\chi_{\Gamma_6\backslash \Omega},
\widetilde{M}_{8}$ that if $|\xi_3^*|\ll N$, then
$$
|M_{6}\cdot\chi_{\Gamma_6\backslash \Omega}|\lesssim {N_1^*}^{\frac{1}{2}} {N_3^*}^{\frac{1}{2}}N_4^*,\quad
|\widetilde{M}_{8}|\lesssim {N_1^*}^{\frac{1}{2}} {N_3^*}^{\frac{1}{2}}
$$
See Corollary 4.1 and Proposition 4.3 below.

\item[\rm (c)] For the construction of $\Omega_3$, we have two
observations. On one hand, we can prove (see Lemma 4.9 below) that
$$
|\alpha_6|\sim {N_1^*}^2, \mbox{ for } (\xi_1,\cdots,\xi_6)\in \Omega_3.
$$
On the other hand, it rules out the bad case
$$
|\xi_3^*|\gtrsim N\gg |\xi_4^*|
$$
in the resonant set $\Gamma_6\backslash \Omega$. This case prevents us to
give a better 6-linear estimate, see Proposition 5.1 below.
\end{enumerate}
\end{remark}

\section{Upper bound of the multipliers: $M_6, M_8, \widetilde{M}_8, M_{10}$}

The key ingredient to prove the almost conservation properties of
the modified energies is to obtain the upper bounds of the
multipliers introduced in Section 3. In this section, we will
present a detailed analysis of the multipliers: $M_6, M_8,
\widetilde{M}_8, M_{10}$.

\subsection{An alternative description of the multipliers: $M_6, M_8, \widetilde{M}_8$}
As a preparation of the next subsections, we rewrite the multipliers in a bright way by merging similar items.

\begin{lem}\label{M6-R}
For the multiplier $M_6$ defined in  (\ref{M6}), we have
$$
M_6=\beta_6+I_1+I_2+I_3+I_4+I_5+I_6,
$$
where  $\beta_6$ defined in (\ref{beta-6}) and
\begin{equation*}\aligned
I_1=&\; C_6\left[M_4(\xi_3,\xi_{214},\xi_5,\xi_6)+M_4(\xi_3,\xi_{216},\xi_5,\xi_4)
+M_4(\xi_3,\xi_{416},\xi_5,\xi_2)\right]\xi_1, \\
I_2=&\; C_6\left[M_4(\xi_{123},\xi_{4},\xi_5,\xi_6)+M_4(\xi_{125},\xi_{4},\xi_3,\xi_6)
+M_4(\xi_{325},\xi_{4},\xi_1,\xi_6)\right]\xi_2, \\
I_3=&\; C_6\left[M_4(\xi_1,\xi_{234},\xi_5,\xi_6)+M_4(\xi_1,\xi_{236},\xi_5,\xi_4)
+M_4(\xi_1,\xi_{436},\xi_5,\xi_2)\right]\xi_3, \\
I_4=&\; C_6\left[M_4(\xi_{143},\xi_{2},\xi_5,\xi_6)+M_4(\xi_{145},\xi_{2},\xi_3,\xi_6)
+M_4(\xi_{345},\xi_{2},\xi_1,\xi_6)\right]\xi_4, \\
I_5=&\; C_6\left[M_4(\xi_1,\xi_{254},\xi_3,\xi_6)+M_4(\xi_1,\xi_{256},\xi_3,\xi_4)
+M_4(\xi_1,\xi_{456},\xi_3,\xi_2)\right]\xi_5, \\
I_6=&\; C_6\left[M_4(\xi_{163},\xi_{2},\xi_5,\xi_4)+M_4(\xi_{165},\xi_{2},\xi_3,\xi_4)
+M_4(\xi_{365},\xi_{2},\xi_1,\xi_4)\right]\xi_6
\endaligned\end{equation*}
for some constant $C_6$.
\end{lem}

For $M_8$, we rewrite it as the following two formulations.
\begin{lem} For the multiplier $M_8$ defined in  (\ref{M8}), we have
\begin{equation}\label{M8-R}
M_8=J_1+J_2+J_3+J_4=J_1'+J_2'+J_3'+J_4',
\end{equation}
where
\begin{align*}\aligned
J_1=&\; 2C_8'\sum\limits_{{\{a,c,e\}=\{3,5,7\}}\atop{\{b,d,f\}=\{4,6,8\}}}
\left[M_4(\xi_{12abc},\xi_d,\xi_e,\xi_f)-M_4(\xi_a,\xi_{12bcd},\xi_e,\xi_f)\right], \\
J_2=&\; C_8'\sum\limits_{{\{a,c,e\}=\{3,5,7\}}\atop{\{b,d,f\}=\{4,6,8\}}}
\left[M_4(\xi_{a2cbe},\xi_d,\xi_1,\xi_f)-M_4(\xi_a,\xi_{b1dcf},\xi_e,\xi_2)\right], \\
J_3=&\; C_8'\sum\limits_{{\{a,c,e\}=\{3,5,7\}}\atop{\{b,d,f\}=\{4,6,8\}}}
\left[M_4(\xi_{1badc},\xi_2,\xi_e,\xi_f)-M_4(\xi_1,\xi_{2abcd},\xi_e,\xi_f)\right], \\
J_4=&\;
2C_8'\sum\limits_{{\{a,c,e\}=\{3,5,7\}}\atop{\{b,d,f\}=\{4,6,8\}}}
\left[M_4(\xi_{1},\xi_2,\xi_{abcde},\xi_f)-M_4(\xi_1,\xi_2,\xi_a,\xi_{bcdef})\right], \\
J_1'=&\;
2C_8'\sum\limits_{{\{a,c\}=\{5,7\}}\atop{\{b,d,f,h\}=\{2,4,6,8\}}}
\left[M_4(\xi_{1b3da},\xi_f,\xi_c,\xi_h)-M_4(\xi_{a},\xi_{b1d3f},\xi_c,\xi_h)\right], \\
J_2'=&\; C_8'\sum\limits_{{\{a,c\}=\{5,7\}}\atop{\{b,d,f,h\}=\{2,4,6,8\}}}
\left[M_4(\xi_{1badc},\xi_f,\xi_3,\xi_h)+M_4(\xi_{3badc},\xi_f,\xi_1,\xi_h)\right], \\
J_3'=&\; -C_8'\sum\limits_{{\{a,c\}=\{5,7\}}\atop{\{b,d,f,h\}=\{2,4,6,8\}}}
\left[M_4(\xi_1,\xi_{b3daf},\xi_c,\xi_f)+M_4(\xi_3,\xi_{b1daf},\xi_c,\xi_f)\right], \\
J_4'=&\;
2C_8'\sum\limits_{{\{a,c,e\}=\{3,5,7\}}\atop{\{b,d,f\}=\{4,6,8\}}}
\left[M_4(\xi_{1},\xi_{badcf},\xi_3,\xi_e)-M_4(\xi_1,\xi_b,\xi_3,\xi_{adcfe})\right] \\
\endaligned\end{align*}
for some constant $C_8'$.
\end{lem}
For $\widetilde{M}_8$, we rewrite it as follows.
\begin{lem} For the multiplier $\widetilde{M}_8$ defined in  (\ref{M8'}), we have
\begin{equation}\label{M8'-R}
\widetilde{M}_8=\tilde{J}_1+\tilde{J}_2+\tilde{J}_3+\tilde{R}_8,
\end{equation}
where
\begin{align*}\aligned
\tilde{J}_1=&\; \tilde{C}_8'
\Big[\sigma_6(\xi_3,\xi_{214},\xi_5,\xi_6,\xi_7,\xi_8)+\sigma_6(\xi_3,\xi_{216},\xi_5,\xi_4,\xi_7,\xi_8)\\
&\;
+\sigma_6(\xi_3,\xi_{218},\xi_5,\xi_4,\xi_7,\xi_6)+\sigma_6(\xi_3,\xi_{416},\xi_5,\xi_2,\xi_7,\xi_8)\\
&\;
+\sigma_6(\xi_3,\xi_{418},\xi_5,\xi_2,\xi_7,\xi_6)+\sigma_6(\xi_3,\xi_{618},\xi_5,\xi_2,\xi_7,\xi_4)
\Big]\xi_1, \\
\tilde{J}_2=&\; \tilde{C}_8'
\Big[\sigma_6(\xi_{123},\xi_4,\xi_5,\xi_6,\xi_7,\xi_8)+\sigma_6(\xi_{125},\xi_4,\xi_3,\xi_6,\xi_7,\xi_8)\\
&\;
+\sigma_6(\xi_{127},\xi_4,\xi_3,\xi_6,\xi_5,\xi_8)+\sigma_6(\xi_{325},\xi_4,\xi_1,\xi_6,\xi_7,\xi_8)\\
&\;
+\sigma_6(\xi_{327},\xi_4,\xi_1,\xi_6,\xi_5,\xi_8)+\sigma_6(\xi_{527},\xi_4,\xi_1,\xi_6,\xi_3,\xi_8)
\Big]\xi_2, \\
\tilde{J}_3=&\; \tilde{C}_8'
\Big[\sigma_6(\xi_1,\xi_{234},\xi_5,\xi_6,\xi_7,\xi_8)+\sigma_6(\xi_1,\xi_{236},\xi_5,\xi_4,\xi_7,\xi_8)\\
&\;
+\sigma_6(\xi_1,\xi_{238},\xi_5,\xi_4,\xi_7,\xi_6)+\sigma_6(\xi_1,\xi_{436},\xi_5,\xi_2,\xi_7,\xi_8)\\
&\;
+\sigma_6(\xi_1,\xi_{438},\xi_5,\xi_2,\xi_7,\xi_6)+\sigma_6(\xi_1,\xi_{638},\xi_5,\xi_2,\xi_7,\xi_4)
\Big]\xi_3
\endaligned\end{align*}
for some constant $\tilde{C}_8'$, and
\begin{equation}\label{R8-tilde}
\big|\tilde{R}_8\big| \lesssim \max_{\Omega} |\sigma_6|\cdot
\max\{|\xi_4|,\cdots, |\xi_8|\} .
\end{equation}
\end{lem}

Next, we give the bounds of the multipliers one by one. First,
we may assume by symmetry that
$$
|\xi_1|\geq |\xi_2|
$$
in the following analysis. Hence $$\xi_1^{*}=\xi_1, \;\;
\xi_2^*=\xi_2 \mbox{ or } \xi_3.$$

\subsection{Known facts}
In this subsection, we restate some results obtained in
\cite{CKSTT-02-DNLS}. First, we have
\begin{lem}[\cite{CKSTT-02-DNLS}] If $N_1^*\ll N$, then we have
\begin{eqnarray}
&M_4(\xi_1,\xi_2,\xi_3,\xi_4)=\frac{1}{2}(\xi_1+\xi_3),\label{4.40}\\
&M_6(\xi_1,\cdots,\xi_6)=0, \quad M_8(\xi_1,\cdots,\xi_8)=0.\label{4.41}
\end{eqnarray}
\end{lem}
Second, we present some estimates on the multipliers.
\begin{lem}[\cite{CKSTT-02-DNLS}]
The following estimates hold:
\begin{itemize}
\item[{\rm(1)}]
           \begin{equation}\label{EM4-0}
           |M_{4}(\xi_1,\xi_2,\xi_3,\xi_4)|\lesssim m_1^2 N_1^*;
           \end{equation}
\item[{\rm(2)}]
           If $|\xi_1|\sim |\xi_3|\gtrsim N\gg |\xi_3^*|$, then
           \begin{equation}\label{EM4-1}
           |M_{4}(\xi_1,\xi_2,\xi_3,\xi_4)|\lesssim m_1^2 N_3^*;
           \end{equation}
\item[{\rm(3)}]
           If $|\xi_1|\sim |\xi_2|\gtrsim N\gg |\xi_3^*|$, then
           \begin{equation}\label{EM4-2}
           M_{4}(\xi_1,\xi_2,\xi_3,\xi_4)= \frac{1}{2}m_1^2 \xi_1+ R(\xi_1,\xi_2,\xi_3,\xi_4), \quad\mbox{ for } |R|\lesssim N_3^*.
           \end{equation}
\item[{\rm(4)}]
            If $|\xi_3^*|\gtrsim N$, then
           \begin{equation}\label{EM6-1}
           |M_{6}(\xi_1,\cdots,\xi_6)|\lesssim m_1^2 {N_1^*}^2;
           \end{equation}
\item[{\rm(5)}]
           If $|\xi_3^*|\ll N$, then
           \begin{equation}\label{EM6-2}
           |M_{6}(\xi_1,\cdots,\xi_6)|\lesssim N_1^* N_3^*.
           \end{equation}
\end{itemize}
\end{lem}

\subsection{An improvement upper bound of $M_6$}
The estimates (\ref{EM6-1}) and (\ref{EM6-2}) are not enough for us
to use, now we make some refinements.
\begin{prop}\label{EM6-fur} For the multiplier $M_6$ defined in (\ref{M6}), the following estimates hold:
\begin{itemize}
\item[{\rm(1)}] If $\xi_2^*=\xi_2$, then
           \begin{equation}\label{EM6-1-fur}
           |M_{6}(\xi_1,\cdots,\xi_6)|\lesssim N_1^*N_3^*.
           \end{equation}
\item[{\rm(2)}]
           If $|\xi_1|\sim |\xi_2|\gtrsim N\gg |\xi_3^*|$, then
           \begin{align}\label{EM6-2-fur}
           M_{6}(\xi_1,\cdots,\xi_6)
            =-C_6\xi_{1}\xi_{12}+C_6' (m_2^2\xi_2^2-m_1^2\xi_1^2)-C_6m_1^2\xi_1\xi_{12}+O({N_3^*}^2),
           \end{align}
           where $C_6$ is the constant in Lemma \ref{M6-R} and $C_6'=\frac{1}{2}C_6-\frac{i}{6}$.
\end{itemize}
\end{prop}
\begin{proof} For the sake of simplicity, we may assume that $C_6=1$. Further, for
(\ref{EM6-1-fur}), we only consider the case $N_1^*\gg N_3^*$,
otherwise it is contained in (\ref{EM6-1}). Thus, we may assume that
$|\xi_1|\sim |\xi_2|\gg |\xi_3^*|$ in (1).

Now we estimate (\ref{EM6-1-fur}) and (\ref{EM6-2-fur}) together.
Note that
$$
\beta_6=-\frac{i}{6} (m_2^2\xi_2^2-m_1^2\xi_1^2)+O({N_3^*}^2).
$$
It suffices to estimate: $I_1,
\cdots, I_6$ by Lemma \ref{M6-R}.

For $I_1, I_2$, by the definitions, we further divide them into
three parts:
$$
I_1:=I_{11}+I_{12}+I_{13}; \quad I_2:=I_{21}+I_{22}+I_{23},
$$
where
\begin{eqnarray*}
& I_{11}:=M_4(\xi_3,\xi_{214},\xi_5,\xi_6)\xi_1,
\quad I_{12}:=M_4(\xi_3,\xi_{216},\xi_5,\xi_4)\xi_1,\\
&  I_{13}:=M_4(\xi_3,\xi_{416},\xi_5,\xi_2)\xi_1, \\
& I_{21}:=M_4(\xi_{123},\xi_{4},\xi_5,\xi_6)\xi_2,
\quad I_{22}:=M_4(\xi_{125},\xi_{4},\xi_3,\xi_6)\xi_2,\\
&  I_{23}:=M_4(\xi_{325},\xi_{4},\xi_1,\xi_6)\xi_2.
\end{eqnarray*}

In order to estimate  $I_1, \ldots, I_6$, it is enough to prove the
following three lemmas.
\begin{lem} If $|\xi_1|\sim |\xi_2|\gg |\xi_3^*|$,  then we have
\begin{equation}\label{I13addI23}
I_{13}+I_{23}= \frac{1}{2}
(m_1^2\xi_1\xi_2+m_2^2\xi_2^2)+O({N_3^*}^2).
\end{equation}
Hence,
\begin{equation}\label{I13addI23-2}
|I_{13}+I_{23}|\lesssim N_1^* N_3^*.
\end{equation}
\end{lem}
\begin{proof} By the definition, we have
\begin{equation*}\aligned
I_{13}&\;=M_4(\xi_3,\xi_{416},\xi_5,\xi_2)\xi_1\\
&\;=-\frac{m_{416}^2\xi_{416}^2\xi_2+m_2^2\xi_2^2\xi_{416}+m_3^2\xi_3^2\xi_5+m_5^2\xi_5^2\xi_3}
    {\alpha}\cdot \xi_1,
\endaligned\end{equation*}
where
$\alpha= \xi_3^2-\xi_{416}^2+\xi_5^2-\xi_2^2$. Similarly,
\begin{equation*}\aligned
I_{23}&\;=M_4(\xi_{325},\xi_4,\xi_1,\xi_6)\xi_2\\
&\;=-\frac{m_{325}^2\xi_{325}^2\xi_1+m_1^2\xi_1^2\xi_{325}+m_4^2\xi_4^2\xi_6+m_6^2\xi_6^2\xi_4}
    {\alpha'}\cdot \xi_2,
\endaligned\end{equation*}
where
$\alpha'= \xi_{325}^2-\xi_{4}^2+\xi_1^2-\xi_6^2$. Note that $|\xi_1|\sim |\xi_2|\gg |\xi_3^*|$, we have
\begin{equation}\label{al-al'}
|\alpha|,|\alpha'|\sim {N_1^*}^2,
\end{equation}
Then,
\begin{equation*}\aligned
I_{13}&\;=-\frac{m_{416}^2\xi_{416}^2\xi_2+m_2^2\xi_2^2\xi_{416}}
    {\alpha}\cdot \xi_1+O({N_3^*}^2N_4^*/N_1^*),\\
I_{23}&\;=-\frac{m_{325}^2\xi_{325}^2\xi_1+m_1^2\xi_1^2\xi_{325}}
    {\alpha'}\cdot \xi_2+O({N_3^*}^2N_4^*/N_1^*),
\endaligned\end{equation*}
which yield that
\begin{align}\label{I13addI23'}
I_{13}+I_{23}=&-\frac{m_{416}^2\xi_{416}^2\xi_2+m_2^2\xi_2^2\xi_{416}}{\alpha}\cdot \left(\xi_1+\xi_2\right)\\
 & + \xi_2\cdot
\Big(\frac{m_{416}^2\xi_{416}^2\xi_2+m_2^2\xi_2^2\xi_{416}}{\alpha}
-\frac{m_{325}^2\xi_{325}^2\xi_1+m_1^2\xi_1^2\xi_{325}}{\alpha'}\Big)\nonumber\\
 &
+O({N_3^*}^2N_4^*/N_1^*)\nonumber\\
 :=& II_1+\xi_2\cdot II +O({N_3^*}^2N_4^*/N_1^*). \nonumber
\end{align}
First, by the mean value theorem (\ref{MTV}) and $m\leq 1$, we have
\begin{equation}\label{II1}
    |II_1|\lesssim m_1^2|\xi_1+\xi_2||\xi_{1246}|\lesssim {N_3^*}^2.
\end{equation}
On the other hand, note that $\xi_{416}=-\xi_{325}$, we have
\begin{align}\label{II-1-2}
II=&\frac{m_{416}^2\xi_{416}^2\xi_2+m_2^2\xi_2^2\xi_{416}}{\alpha}
-\frac{m_{325}^2\xi_{325}^2\xi_1+m_1^2\xi_1^2\xi_{325}}{\alpha'}\\
 =&
\frac{1}{\alpha}\left(m_{416}^2\xi_{416}^2\xi_2+m_2^2\xi_2^2\xi_{416}
+m_{325}^2\xi_{325}^2\xi_1+m_1^2\xi_1^2\xi_{325}\right)\nonumber\\
 &
-\left(\frac{\alpha+\alpha'}{\alpha\alpha'}\right)\cdot\left(m_{325}^2\xi_{325}^2\xi_1+m_1^2\xi_1^2\xi_{325}\right)\nonumber\\
 =& \frac{1}{\alpha}\Big(m_{416}^2\xi_{416}^2(\xi_1+\xi_2)
+\xi_{416}\left(m_2^2\xi_2^2-m_1^2\xi_1^2\right)\Big)\nonumber\\
 &
-\left(\frac{\alpha+\alpha'}{\alpha\alpha'}\right)\cdot\left(m_{325}^2\xi_{325}^2\xi_1+m_1^2\xi_1^2\xi_{325}\right)
.\nonumber
\end{align}
By the mean value theorem (\ref{MTV}), we have
$$
\frac{1}{\alpha}= \frac{1}{2\xi_1\xi_2}+O(N_3^*/{N_1^*}^3),
$$
$$
m_{416}^2\xi_{416}^2(\xi_1+\xi_2)
+\xi_{416}\left(m_2^2\xi_2^2-m_1^2\xi_1^2\right)=O({N_1^*}^2 N_3^*).
$$
Thus,
\begin{align} \label{II-1}
\mbox{Term 1 of } & (\ref{II-1-2}) \\
=&\; \frac{1}{2\xi_1\xi_2}\Big(m_{416}^2\xi_{416}^2(\xi_1+\xi_2)
+\xi_{416}\left(m_2^2\xi_2^2-m_1^2\xi_1^2\right)\Big)
+O({N_3^*}^2/N_1^*)\nonumber\\
 =&\; \frac{1}{2\xi_1\xi_2}\Big(m_{1}^2\xi_{1}^2(\xi_1+\xi_2)
+\xi_{1}\left(m_2^2\xi_2^2-m_1^2\xi_1^2\right)\Big)+O({N_3^*}^2/N_1^*)\nonumber\\
 =&\; \frac{1}{2}\left(m_{1}^2\xi_{1}
+m_2^2\xi_2\right)+O({N_3^*}^2/N_1^*) .\nonumber
\end{align}
On the other hand, by mean value theorem (\ref{MTV}),
$$
|\alpha+\alpha'|=|\alpha_6|\lesssim N_1^* N_3^*,\quad |m_{325}^2\xi_{325}^2\xi_1+m_1^2\xi_1^2\xi_{325}|=O({N_1^*}^2 N_3^*).
$$
Thus, by (\ref{al-al'}), we get
\begin{equation}
\mbox{Term 2 of (\ref{II-1-2})}\lesssim  {N_3^*}^2 /{N_1^*}. \label{II-2}
\end{equation}
Combining (\ref{II-1-2}), (\ref{II-1}) with (\ref{II-2}), we have
\begin{equation}\label{E-II}
II= \frac{1}{2}\left(m_{1}^2\xi_{1}
+m_2^2\xi_2\right)+O({N_3^*}^2/N_1^*).
\end{equation}
Inserting (\ref{II1}) and (\ref{E-II}) into (\ref{I13addI23'}), we
have the desired result.
\end{proof}

\begin{lem} If $|\xi_1|\sim |\xi_2|\gg |\xi_3^*|$,  then we have
\begin{equation}\label{I11addI22}
I_{11}+I_{12}+I_{21}+I_{22}\lesssim N_1^* N_3^*.
\end{equation}
Furthermore, if $|\xi_3^*|\ll N$,  we have
\begin{equation}\label{I11addI22'}
I_{11}+I_{12}+I_{21}+I_{22}=-\xi_1\xi_{12}+O({N_3^*}^2).
\end{equation}
\end{lem}
\begin{proof} Since $|\xi_{12}|\lesssim N_3^*,$ (\ref{I11addI22}) follows from (\ref{EM4-0}). Moreover,
if $|\xi_{12}|\ll N$, then by (\ref{4.40}), we have
$$
I_{11}=I_{12}=\frac{1}{2}\xi_{35}\cdot\xi_1, \quad
I_{21}=I_{22}=-\frac{1}{2}\xi_{46}\cdot\xi_2,
$$
which imply that
$$
I_{11}+I_{12}+I_{21}+I_{22}=\xi_{35}\xi_1-\xi_{46}\xi_2=\xi_{12}\cdot\xi_{235}=-\xi_1\xi_{12}+O({N_3^*}^2).
$$
This completes the proof of the lemma.
\end{proof}

\begin{lem} If $|\xi_1|\sim |\xi_2|\gg |\xi_3^*|$, then we have
\begin{equation}\label{I3addI6'}
|I_{3}+I_{4}+I_{5}+I_{6}|\lesssim N_1^* N_3^*.
\end{equation}
Furthermore, if $|\xi_3^*|\ll N$, then
\begin{equation}\label{I3addI6}
I_{3}+I_{4}+I_{5}+I_{6}=
-\frac{3}{2}m_1^2\xi_1\xi_{12}+O({N_3^*}^2).
\end{equation}
\end{lem}
\begin{proof}
(\ref{I3addI6'}) follows from (\ref{EM4-0}). Now we
consider the case $|\xi_3^*|\ll N$. By (\ref{EM4-2}), we have
\begin{equation}\label{EM4-3}
M_{4}(\overline{\xi_1},\overline{\xi_2},\xi_3,\xi_4)= \frac{1}{2}m_1^2 \xi_1+ O(N_3^*),
\end{equation}
where $ \overline{\xi_1}+\overline{\xi_2}+\xi_3+\xi_4=0,\;
\overline{\xi_1}=\xi_1+O(N_3^*),\; \overline{\xi_2}=\xi_2+O(N_3^*) $
and $|\overline{\xi_1}|\sim |\overline{\xi_2}|\gtrsim N \gg |\xi_3^*|$. Using
(\ref{EM4-3}), we obtain
\begin{equation*}\aligned
I_{3}+I_{4}+I_{5}+I_{6} &\;= \frac{3}{2}m_1^2\xi_1(\xi_3+\xi_4+\xi_5+\xi_6)+O({N_3^*}^2)\\
 &\;= -\frac{3}{2}m_1^2\xi_1(\xi_1+\xi_2)+O({N_3^*}^2).
\endaligned\end{equation*}
This completes the proof of the lemma.
\end{proof}

Now we finish the proof of Proposition \ref{EM6-fur}.  Indeed,
(\ref{EM6-1-fur}) follows from (\ref{I13addI23-2}),
(\ref{I11addI22}) and (\ref{I3addI6'}). While by (\ref{I13addI23})
and (\ref{I3addI6}), we have,
\begin{align}\label{I13-I6}
I_{13}+I_{23}+I_{3} & +I_{4}+I_{5}+I_{6} \\
 &= \frac{1}{2}
(m_1^2\xi_1\xi_2+m_2^2\xi_2^2)-\frac{3}{2}m_1^2\xi_1(\xi_1+\xi_2)
+O({N_3^*}^2)\nonumber\\
 &= \frac{1}{2} (m_2^2\xi_2^2-m_1^2\xi_1^2) -m_1^2\xi_1(\xi_1+\xi_2)
+O({N_3^*}^2).\nonumber
\end{align}
Therefore, (\ref{EM6-2-fur}) follows from
(\ref{I11addI22'}) and (\ref{I13-I6}).
\end{proof}
\begin{cor} If $ |\xi_3^*|\ll N$, then we have
\begin{equation}\label{EM6-3}
    |M_{6}(\xi_1,\cdots,\xi_6)|\lesssim {N_1^*}^{\frac{1}{2}} {N_3^*}^{\frac{1}{2}}N_4^* \;\; \text{in}\,\,  \Gamma_6\backslash \Omega.
\end{equation}
\end{cor}
\begin{proof} In this situation, $\xi_2^*=\xi_2$ (see Remark 3.3 (a)). Then by (\ref{EM6-2-fur}) and the mean value theorem (\ref{MTV}), we
have
$$
|M_{6}(\xi_1,\cdots,\xi_6)|\lesssim |\xi_1||\xi_1+\xi_2|+{N_3^*}^2.
$$
Moreover,  since $|\xi_1|^{\frac{1}{2}}|\xi_1+\xi_2|\lesssim
|\xi_3^*|^{\frac{3}{2}}$ in $ \Gamma_6\backslash\Omega$, we have
$$
|M_{6}(\xi_1,\cdots,\xi_6)|\lesssim {N_1^*}^{\frac{1}{2}} {N_3^*}^{\frac{3}{2}}.
$$
Then (\ref{EM6-3}) follows by the fact that $N_3^*\sim N_4^*$ in $
\Gamma_6\backslash \Omega_3$.
\end{proof}

\subsection{A upper bound of $M_8$}
\begin{prop}\label{EM8-Prop}
\begin{equation}\label{EM8-1}
|M_8(\xi_1,\cdots,\xi_8)|\lesssim N_1^*.
\end{equation}
Furthermore, if $|\xi_3^*|\ll N$, then we have
\begin{equation}\label{EM8-2}
|M_8(\xi_1,\cdots,\xi_8)|\lesssim N_3^*.
\end{equation}
\end{prop}
\begin{proof} By (\ref{EM4-0}), we have $|M_4(\xi_1,\xi_2,\xi_3,\xi_4)|\lesssim N_1^*$. Thus (\ref{EM8-1}) follows.
For (\ref{EM8-2}), we split it into two cases.

\noindent\textbf{Case 1, $\xi_2^*=\xi_2$. } By (\ref{M8-R}), we have
$$
M_8=J_1+J_2+J_3+J_4.
$$
So it suffices to prove: $|J_1|,|J_2|,|J_3|,|J_4|\lesssim N_3^*.$
First, $J_1$ follows immediately from $|\xi_1+\xi_2|\lesssim N_3^*$
and (\ref{EM4-0}). While $J_2$ follows from (\ref{EM4-1}) and $J_3,
J_4$ follow from (\ref{EM4-3}).

\noindent\textbf{Case 2, $\xi_2^*=\xi_3$.} Now we adopt the
formulation:
$$
M_8=J_1'+J_2'+J_3'+J_4',
$$
and it is necessary to prove: $|J_1'|,|J_2'|,|J_3'|,|J_4'|\lesssim
N_3^*.$  $J_1'$ and $ J_2'$ are similar to $J_1$ and $J_2$. For
$J'_3$, we also use (\ref{EM4-3}) to give
$$
J_3'= C(m_1^2\xi_1+m_3^2\xi_3)+O(N_3^*)= O(N_3^*),
$$
where we used the mean value theorem (\ref{MTV}). $J_4'$ is similar
to $J_2$.
\end{proof}

\subsection{A upper bound of $\sigma_6, \widetilde{M}_8$}
First, we prove that $\sigma_6$ is uniformly bounded in $\Omega$, which implies that the set $\Omega$
is non-resonant.
\begin{lem}\label{Sigma-Lem} In $\Omega$, we have
           \begin{equation}\label{sigma6-1}
           |\sigma_{6}(\xi_1,\cdots,\xi_6)|\lesssim 1.
           \end{equation}
           Particularly, in $\Omega_1\cap \big\{|\xi_3^*|\ll
           N\big\}$, we have
           \begin{equation}\label{sigma6-2}
           |\sigma_{6}(\xi_1,\cdots,\xi_6)|\lesssim N_3^*/N_1^*.
           \end{equation}
\end{lem}
\begin{proof}
Recall that
$$
\sigma_6=-\frac{M_6}{\alpha_6}\cdot \chi_{\Omega}, \quad \alpha_6=-i(\xi_1^2-\xi_2^2+\xi_3^2-\xi_4^2+\xi_5^2-\xi_6^2).
$$

In $\Omega_1$, we have
$$
|\alpha_6(\xi_1,\cdots,\xi_6)|\sim {N_1^*}^2.
$$
This gives (\ref{sigma6-2}) by (\ref{EM6-2}) and (\ref{sigma6-1}) by (\ref{EM6-1}).

In $\Omega_2$, we have
$$
|\xi_1^2-\xi_2^2|\sim |\xi_1||\xi_1+\xi_2|\gg |\xi_3^*|^2,
$$
which yields that
\begin{equation}\label{alfa-6}
|\alpha_6|\sim |\xi_1||\xi_1+\xi_2|.
\end{equation}
While from (\ref{EM6-2-fur}) and the mean value theorem (\ref{MTV}), we have
$$
|M_{6}(\xi_1,\cdots,\xi_6)|\lesssim |\xi_1||\xi_1+\xi_2|+ {N_3^*}^2\lesssim |\xi_1||\xi_1+\xi_2|.
$$
This gives (\ref{sigma6-1}) in $\Omega_2$.

In $\Omega_3$, since $\xi_1^*\cdot \xi_2^*<0, \xi_2^*\cdot
\xi_3^*>0$, it holds that
$$
|\xi_1^*|= |\xi_2^*|+ |\xi_3^*|+o(N_3^*).
$$
We claim that
\begin{equation}\label{L4.9-1}
|\alpha_6|\gtrsim N_1^* N_3^*.
\end{equation}
Indeed, for (\ref{L4.9-1}), we divide into the following three
cases:
$$
\mbox{(i)}\, \xi_2^*=\xi_2, \xi_3^*=\xi_3;\quad \mbox{(ii)}\,
\xi_2^*=\xi_2, \xi_3^*=\xi_4;\quad \mbox{(iii)}\,
\xi_2^*=\xi_3, \xi_3^*=\xi_2.
$$

If $\xi_2^*=\xi_2, \xi_3^*=\xi_3$, then we get
\begin{equation*}\aligned
|\alpha_6|=& \;\big|(\xi_1^2-\xi_2^2)+\xi_3^2+(-\xi_4^2+\xi_5^2-\xi_6^2)\big|\\
  =& \;\big(\xi_1^2-\xi_2^2\big)+\xi_3^2+o(|\xi_{3}|^2)\\
  =& \;-\xi_1\xi_3+\xi_3^2+o(|\xi_{1}||\xi_{3}|)\\
  \sim& \;|\xi_1||\xi_3|.
\endaligned\end{equation*}

If $\xi_2^*=\xi_2, \xi_3^*=\xi_4$, then we have
\begin{equation*}\aligned
|\alpha_6|=& \;\big|(\xi_1^2-\xi_2^2-\xi_4^2)+(\xi_3^2+\xi_5^2-\xi_6^2)\big|\\
  =& \;\big(\xi_1^2-\xi_2^2-\xi_4^2\big)+o(|\xi_{4}|^2)\\
  =& \;\left(\big[|\xi_2|+|\xi_4|+o(|\xi_4|)\big]^2-\xi_2^2-\xi_4^2\right)+o(|\xi_{4}|^2)\\
  \sim & \; |\xi_2||\xi_4|.
\endaligned\end{equation*}

If $\xi_2^*=\xi_3, \xi_3^*=\xi_2$, then we have
$$
|\alpha_6|=\big(\xi_1^2-\xi_2^2+\xi_3^2\big)+o(|\xi_{3}|^2)\geq
\xi_3^2+o(|\xi_{3}|^2)\sim \xi_1^2.
$$
This proves (\ref{L4.9-1}).

By (\ref{EM6-1}) and (\ref{EM6-2}), we have
$|M_{6}(\xi_1,\cdots,\xi_6)|\lesssim {N_1^*}^2$. Then
(\ref{sigma6-1}) follows if $N_1^*\sim N_3^*$. Now we consider the
other case: $N_1^*\gg N_3^*$. Thus we have: $\xi_2^*=\xi_2$ in
$\Omega_3\backslash\Omega_1$. Then  (\ref{sigma6-1}) in $\Omega_3\backslash
\Omega_1$ follows from (\ref{EM6-1-fur}) and (\ref{L4.9-1}).
\end{proof}

Now we give the upper bound of $\widetilde{M}_8$.
\begin{prop}\label{EM8'-Prop}
\begin{equation}\label{EM8'-1}
|\widetilde{M}_8(\xi_1,\cdots,\xi_8)|\lesssim N_1^*.
\end{equation}
Furthermore, if $|\xi_3^*|\ll N$, then we have
\begin{equation}\label{EM8'-2}
|\widetilde{M}_8(\xi_1,\cdots,\xi_8)|\lesssim {N_1^*}^{\frac{1}{2}}{N_3^*}^{\frac{1}{2}}.
\end{equation}
\end{prop}
\begin{proof}
Since $|\sigma_{6}|\lesssim 1,$ we have (\ref{EM8'-1}). Now we
turn to (\ref{EM8'-2}).  By (\ref{M8'-R}), we
shall estimates: $\tilde{J}_1, \tilde{J}_2,\tilde{J}_3$. For this purpose,
we divide it into two cases.

\noindent\textbf{Case 1, $\xi_2^*=\xi_2$.} Since
$|\sigma_{6}|\lesssim 1,$ we have $|\tilde{J}_3|\lesssim N_3^*$. Now
we consider the other two parts. Since $\sigma_6=0$ for
$|\xi_1^*|\ll N$, we know that the first, second, third terms of
$\tilde{J}_1, \tilde{J}_2$ vanish. Therefore,
\begin{align}\label{M8'-2}
\widetilde{M}_8 =&
\tilde{C}_8'\big[\sigma_6(\xi_3,\xi_{416},\xi_5,\xi_2,\xi_7,\xi_8)
+\sigma_6(\xi_3,\xi_{418},\xi_5,\xi_2,\xi_7,\xi_6)\\
 &+\sigma_6(\xi_3,\xi_{618},\xi_5,\xi_2,\xi_7,\xi_4)
\big]\xi_1+\tilde{C}_8'\big[\sigma_6(\xi_{325},\xi_4,\xi_1,\xi_6,\xi_7,\xi_8) \nonumber\\
 &
+\sigma_6(\xi_{327},\xi_4,\xi_1,\xi_6,\xi_5,\xi_8)+\sigma_6(\xi_{527},\xi_4,\xi_1,\xi_6,\xi_3,\xi_8)
\big]\xi_2+O(N_3^*).\nonumber
\end{align}
By (\ref{sigma6-2}), each term is bounded by $N_3^*$.

\noindent \textbf{Case 2, $\xi_2^*=\xi_3$.} In this case,
$|\tilde{J}_2|\lesssim N_3^*$, so we only need to estimate
$\tilde{J}_1, \tilde{J}_3$. By permutating the terms in
$\tilde{J}_1, \tilde{J}_3$, we may rewrite $\widetilde{M}_8$ as
\begin{equation*}
\aligned
\widetilde{M}_8=&\;\sum\limits_{{\{a,c\}=\{5,7\}}\atop{\{b,d,f,h\}=\{2,4,6,8\}}}
\big[\sigma_6(\xi_3,\xi_{b1d},\xi_a,\xi_f,\xi_c,\xi_h)\xi_1+
\sigma_6(\xi_1,\xi_{b3d},\xi_a,\xi_f,\xi_c,\xi_h)\xi_3\big]\\
&\;\qquad\qquad\qquad\quad+O(N_3^*).
\endaligned
\end{equation*}
As an example, we only consider
$$
\sigma_6(\xi_3,\xi_{214},\xi_5,\xi_6,\xi_7,\xi_8)\xi_1+
\sigma_6(\xi_1,\xi_{234},\xi_5,\xi_6,\xi_7,\xi_8)\xi_3,
$$
which equals to
\begin{equation}\label{4.72}
    III\cdot \xi_1+ O(N_3^*),
\end{equation}
where
$$
III:=\sigma_6(\xi_3,\xi_{214},\xi_5,\xi_6,\xi_7,\xi_8)-
\sigma_6(\xi_1,\xi_{234},\xi_5,\xi_6,\xi_7,\xi_8).
$$
We first adopt some notations for short. We denote
\begin{equation*}
\aligned
A:=M_6(\xi_3,\xi_{214},\xi_5,\xi_6,\xi_7,\xi_8);&\;\quad A':=M_6(\xi_1,\xi_{234},\xi_5,\xi_6,\xi_7,\xi_8),\\
B:=\alpha_6(\xi_3,\xi_{214},\xi_5,\xi_6,\xi_7,\xi_8);&\; \quad B':=\alpha_6(\xi_1,\xi_{234},\xi_5,\xi_6,\xi_7,\xi_8).
\endaligned
\end{equation*}
Since
$$
\Omega_2(\xi_3,\xi_{214},\xi_5,\xi_6,\xi_7,\xi_8)=\Omega_2(\xi_1,\xi_{234},\xi_5,\xi_6,\xi_7,\xi_8),
$$
then by (\ref{sigma6-1}), (\ref{alfa-6}) and the definition of $\Omega_2$, we have
\begin{equation}\label{4.74}
\left|\frac{A}{B}\right|, \left|\frac{A'}{B'}\right|\lesssim 1; \quad |B|,|B'|\sim |\xi_{1234}||\xi_1|
\gg {N_1^*}^{\frac{1}{2}}{N_3^*}^{\frac{3}{2}}.
\end{equation}
Moreover,
\begin{equation}\label{4.75}
III=\frac{A}{B}-\frac{A'}{B'}=\frac{1}{B}(A+A')-\frac{A'}{B'}\cdot \frac{B+B'}{B}.
\end{equation}
On one hand, by (\ref{EM6-2-fur}) and  (\ref{4.74}), we have
\begin{equation*}
\aligned
A+A'=&\;C_6\xi_{1234}\cdot(2\xi_{2457}+\xi_{13})-C_6\xi_{1234}(m_1^2\xi_1+m_3^2\xi_3)
\\
&\; +C_6'
(m_{214}^2\xi_{214}^2-m_3^2\xi_3^2+m_{234}^2\xi_{234}^2-m_1^2\xi_1^2)+O({N_3^*}^2).
\endaligned
\end{equation*}
Further, by the mean value theorem (\ref{MTV}) in the second term
and by the double mean value theorem (\ref{DMTV}) in the third term,
we have
\begin{equation}\label{AA'}
|A+A'|\lesssim m_1^2|\xi_{1234}||\xi_{24}|+{N_3^*}^2.
\end{equation}
Therefore, by (\ref{4.74}) and (\ref{AA'}), we have
\begin{align}\label{III1}
\left|\frac{1}{B}(A+A')\right| \lesssim&
m_1^2\frac{|\xi_{24}|}{|\xi_1|}+\frac{{N_3^*}^2}{{N_1^*}^{\frac{1}{2}}{N_3^*}^{\frac{3}{2}}}\\
 \lesssim&
N_3^*/N_1^*+{N_3^*}^{\frac{1}{2}}/{N_1^*}^{\frac{1}{2}}\lesssim
{N_3^*}^{\frac{1}{2}}/{N_1^*}^{\frac{1}{2}}.\nonumber
\end{align}
On the other hand,
\begin{equation*}
\aligned
|B+B'|=&\;|\xi_1^2-\xi_{234}^2+\xi_3^2-\xi_{214}^2|+O({N_3^*}^2) =
\; 2|\xi_{1234}||\xi_{24}|+O({N_3^*}^2).
\endaligned
\end{equation*}
Therefore, by the similar estimates as those in (\ref{4.74}) and (\ref{III1}), we have
\begin{equation}\label{III2}
    \left|\frac{A'}{B'}\cdot \frac{B+B'}{B}\right|\lesssim
{N_3^*}^{\frac{1}{2}}/{N_1^*}^{\frac{1}{2}}.
\end{equation}
Inserting (\ref{III1}) and (\ref{III2}) into (\ref{4.75}), we have
$$
|III|\lesssim {N_3^*}^{\frac{1}{2}}/{N_1^*}^{\frac{1}{2}}.
$$
which together with (\ref{4.72}) yields (\ref{EM8'-2}).
\end{proof}

\section{An upper bound on the increment of $E_I^3(u(t))$}

By the multilinear correction analysis, the almost conservation law
of $E_I^3(u(t))$ is the key ingredient to establish the global
well-posedness below the energy space. This is made up of the
following 6-linear, 8-linear and 10-linear estimates.

\begin{prop} \label{6-linear}
For any $s\geq\frac{1}{2}$, we have
\begin{equation}
 \left|\displaystyle\int_0^\delta\Lambda_6(M_6\cdot \chi_{\Gamma_6\backslash\Omega};w(t))\,dt\right|
 \lesssim
 N^{-\frac{5}{2}+}\>
 \|Iw\|_{Y_1}^6.
 \label{Lambda6}
\end{equation}
\end{prop}
\begin{proof}
By (\ref{4.41}), when $|\xi_1|,\cdots, |\xi_6|\ll N$, we have
$M_6=0$. Therefore, we may assume that $|\xi_1^*|\sim
|\xi_2^*|\gtrsim N$. Note that
$$
\|\chi_{[0,\delta]}(t)f\|_{X_{0,\frac{1}{2}-}}\lesssim \|f\|_{X_{0,\frac{1}{2}}}
$$
(see Lemma 2.2 in \cite{Li-Wu} for example),  (\ref{Lambda6}) is
reduce to
$$
  \left|\displaystyle\int\Lambda_6(M_6\cdot \chi_{\Gamma_6\backslash \Omega};w(t))\,dt\right|
 \lesssim
 N^{-\frac{5}{2}+}\>
 \|Iw\|_{X_{1,\frac{1}{2}-}}\|Iw\|_{Y_1}^{5}.
$$
But the $0+$ loss is not essential by (\ref{2.18})--(\ref{2.20}) and
(\ref{XE2}) for $q<6$, thus it will not be mentioned. By
Plancherel's identity and
$\widehat{\bar{f}}(\xi,\tau)=\bar{\hat{f}}(-\xi,-\tau)$, we only
need to show that for any $f_j\in Y_0^+,$   $ j=1,3,5 $ and $f_j\in
Y^-_0$,   $j=2,4,6$,
\begin{align}\label{P5.1-1}
& \int_{\Gamma_6\times \Gamma_6} \frac{M_6\cdot
\chi_{\Gamma_6\backslash \Omega}(\xi_1,\cdots,
\xi_6)\widehat{f_1}(\xi_1,\tau_1)\cdots \widehat{f_6}(\xi_6,\tau_6)}
{\langle\xi_1\rangle m(\xi_1)\cdots \langle\xi_6\rangle m(\xi_6)}\\
&\qquad\qquad\qquad \lesssim
N^{-\frac{5}{2}+}\>\|f_1\|_{Y_0^+}\|f_2\|_{Y_0^-}\cdots
\|f_5\|_{Y_0^+}\|f_6\|_{Y_0^-},\nonumber
\end{align}
where $\Gamma_6\times
\Gamma_6=\{(\xi,\tau):\xi_1+\cdots+\xi_6=0,
\tau_1+\cdots+\tau_6=0\}$, $\xi=(\xi_1,\cdots,\xi_6)$,
$\tau=(\tau_1,\cdots,\tau_6)$. Now we divide it into four
regions:
\begin{equation*}\aligned
  A_1 =& \{(\xi,\tau)\in (\Gamma_6\backslash\Omega)\times \Gamma_6:
  |\xi_2^*|\gtrsim N\gg |\xi_3^*|\}, \\
  A_2 =& \{(\xi,\tau)\in (\Gamma_6\backslash\Omega)\times \Gamma_6:
  |\xi_3^*|\gtrsim N\gg |\xi_4^*|\}, \\
  A_3 =& \{(\xi,\tau)\in (\Gamma_6\backslash\Omega)\times \Gamma_6:
  |\xi_4^*|\gtrsim N\gg |\xi_5^*|\}, \\
  A_4 =& \{(\xi,\tau)\in (\Gamma_6\backslash\Omega)\times \Gamma_6:
  |\xi_5^*|\gtrsim N\}.
\endaligned\end{equation*}
In the following, we adopt the notation $f_j^*$ to be one of $f_j$
for $j=1,\cdots, 6$ and satisfy
$\widehat{f_j^*}=\widehat{f_j^*}(\xi_j^*,\tau_j)$.

\noindent\textbf{Estimate in $A_1$.} By the definition of $\Omega$
and (\ref{EM6-3}), in $(\Gamma_6\backslash\Omega)\times \Gamma_6$, we have
$$
|\xi_1|\sim|\xi_2|\gtrsim N\gg|\xi_3^*|, \mbox{\quad and \quad}
|M_6\cdot \chi_{\Gamma_6\backslash \Omega}|\lesssim  {N_1^*}^{\frac{1}{2}}
{N_3^*}^{\frac{1}{2}}N_4^*.
$$
Therefore, by (\ref{2.18})--(\ref{2.20}), we have
\begin{equation*}\aligned
\mbox{ LHS of (\ref{P5.1-1})}\lesssim &\;
 N^{2s-2} \int_{A_1}
\frac{\widehat{f_1}(\xi_1,\tau_1)\cdots \widehat{f_6}(\xi_6,\tau_6)}
{|\xi_1^*|^{2s-\frac{1}{2}} \langle\xi_3^*\rangle^{\frac{1}{2}}\langle\xi_5^*\rangle\langle\xi_6^*\rangle}\\
= &\;
 N^{2s-2} \int_{A_1}
|\xi_1^*|^{-2s-\frac{1}{2}+}\langle\xi_3^*\rangle^{-\frac{1}{2}}\cdot
(|\xi_1^*|^{\frac{1}{2}-}\widehat{f_1^*}\widehat{f_3^*})\>(|\xi_2^*|^{\frac{1}{2}-}\widehat{f_2^*}\widehat{f_4^*})\\
&\;
\qquad\qquad\,\,\cdot(\langle\xi_5^*\rangle^{-1} \widehat{f_5^*} ) (\langle\xi_6^*\rangle^{-1} \widehat{f_6^*} )
\\
\lesssim& \;
N^{-\frac{5}{2}+}\left\|I^{\frac{1}{2}-}_{\pm}(f_1^*,f_3^*)\right\|_{L^2_{xt}}
\left\|I^{\frac{1}{2}-}_{\pm}(f_2^*,f_4^*)\right\|_{L^2_{xt}}\\
&\;
\qquad\,\cdot
\left\|J_x^{-1}f_5^*\right\|_{L^\infty_{xt}}
\left\|J_x^{-1}f_6^*\right\|_{L^\infty_{xt}}\\
\lesssim&\; N^{-\frac{5}{2}+}\> \|f_1\|_{Y_0^+}\|f_2\|_{Y_0^-}\cdots
\|f_5\|_{Y_0^+}\|f_6\|_{Y_0^-},
\endaligned\end{equation*}
where we use the relations that $|\xi_1^*\pm\xi_3^*|\sim |\xi_1^*|$
and $|\xi_2^*\pm\xi_4^*|\sim |\xi_1^*|$.

\noindent\textbf{Estimate in $A_2$.} Note that $A_2=\emptyset$ in
$(\Gamma_6\backslash\Omega_3)\times \Gamma_6$, thus $M_6\cdot \chi_{\Gamma_6\backslash
\Omega}=0$.

\noindent\textbf{Estimate in $A_3$.} By (\ref{EM6-1}), we have
\begin{equation}\label{EM6-1'}
|M_6\cdot \chi_{\Gamma_6\backslash \Omega}|\lesssim m_1^2{N_1^*}^2.
\end{equation}
Therefore, by (\ref{2.18})--(\ref{2.20}) and (\ref{XE5}), we have
\begin{equation*}\aligned
\mbox{ LHS of (\ref{P5.1-1})}\lesssim &\; N^{2s-2}\displaystyle\int_{A_3}
\frac{\widehat{f_1}(\xi_1,\tau_1)\cdots \widehat{f_6}(\xi_6,\tau_6)}
{|\xi_3^*|^{s} |\xi_4^*|^{s}\langle\xi_5^*\rangle\langle\xi_6\rangle}\\
= &\;
 N^{2s-2} \int_{A_1}
|\xi_1^*|^{-\frac{1}{2}+}|\xi_3^*|^{-s} |\xi_4^*|^{-s}\langle\xi_5^*\rangle^{-1}\cdot
(|\xi_1^*|^{\frac{1}{2}-}\widehat{f_1^*}\widehat{f_5^*})\>(|\xi_2^*|^{0-}\widehat{f_2^*})\\
&\;
\qquad\qquad\,\,\cdot(|\xi_3^*|^{0-}\widehat{f_3^*})(|\xi_4^*|^{0-}\widehat{f_4^*})
(\langle\xi_6^*\rangle^{-1} \widehat{f_6^*} )
\\
\lesssim&\;
N^{-\frac{5}{2}+}\left\|I^{\frac{1}{2}-}_{\pm}(f_1^*,f_5^*)\right\|_{L^2_{xt}}
\left\|J_x^{0-}f_2^*\right\|_{L^6_{xt}}
\left\|J_x^{0-}f_3^*\right\|_{L^6_{xt}}\\
&\;
\qquad\,\cdot
\left\|J_x^{0-}f_4^*\right\|_{L^6_{xt}}
\left\|J_x^{-1}f_6^*\right\|_{L^\infty_{xt}}\\
\lesssim&\; N^{-\frac{5}{2}+}\> \|f_1\|_{Y_0^+}\|f_2\|_{Y_0^-}\cdots
\|f_5\|_{Y_0^+}\|f_6\|_{Y_0^-},
\endaligned\end{equation*}
where we use
the fact that $|\xi_1^*\pm\xi_5^*|\sim
|\xi_1^*|$ in this case.

\noindent\textbf{Estimate in $A_4$.} The worst case is
$|\xi_j|\gtrsim N$ for any $j=1,\cdots, 6$, we only consider this
case. Then by (\ref{EM6-1'}),  (\ref{XE2}) for $q=6-$ and
(\ref{XE6}) for $q=6+$, we have
\begin{equation*}\aligned
\mbox{ LHS of (\ref{P5.1-1})}\lesssim&\; N^{4s-4}\displaystyle\int_{A_4}
\frac{\widehat{f_1}(\xi_1,\tau_1)\cdots \widehat{f_6}(\xi_6,\tau_6)}
{|\xi_3^*|^{s} |\xi_4^*|^{s}|\xi_5^*|^{s}|\xi_6^*|^{s}}\\
\lesssim&\;
N^{-4+}\left\|f_1^*\right\|_{L^{6-}_{xt}}\cdots
\left\|f_5^*\right\|_{L^{6-}_{xt}}
\left\|J_x^{0-}f_6^*\right\|_{L^{6+}_{xt}}\\
\lesssim&\; N^{-4+}\> \|f_1\|_{Y_0^+}\|f_2\|_{Y_0^-}\cdots
\|f_5\|_{Y_0^+}\|f_6\|_{Y_0^-}.
\endaligned\end{equation*}
This gives
the proof of the proposition.
\end{proof}

\begin{prop} \label{8-linear}
For any $s\geq\frac{1}{2}$, we have
\begin{equation}
 \left|\displaystyle\int_0^\delta\Lambda_8(M_8+\widetilde{M}_8;w(t))\,dt\right|
 \lesssim
 N^{-\frac{5}{2}+}\>
 \|Iw\|_{Y_1}^8.
 \label{Lambda8}
\end{equation}
\end{prop}
\begin{proof}
When $|\xi_1|,\cdots, |\xi_8|\ll N$, we have $M_8,\widetilde{M}_8=0$.
Similar to (\ref{P5.1-1}), it suffices to show
\begin{align}\label{P5.2-1}
&\int_{\Gamma_{8}\times \Gamma_{8}}\!\!
\frac{(M_8+\widetilde{M}_8)(\xi_1,\cdots,
\xi_{8})\widehat{f_1}(\xi_1,\tau_1) \cdots
\widehat{f_{8}}(\xi_{8},\tau_{8})} {\langle\xi_1\rangle
m(\xi_1)\cdots \langle\xi_{8}\rangle m(\xi_{8})}\\
&\qquad\qquad\qquad \lesssim
N^{-\frac{5}{2}+}\>\|f_1\|_{Y_0^+}\|f_2\|_{Y_0^-}\cdots
\|f_7\|_{Y_0^+}\|f_8\|_{Y_0^-},\nonumber
\end{align}
where $\Gamma_{8}\times \Gamma_{8}
=\{(\xi_1,\cdots,\xi_{8},\tau_1,\cdots,\tau_{8}):\xi_1+\cdots+\xi_{8}=0,
\tau_1+\cdots+\tau_{8}=0\}$. Now we divide it into three regions:
\begin{equation*}\aligned
  B_1 =& \{(\xi_1,\cdots,\xi_{8},\tau_1,\cdots,\tau_{8})\in \Gamma_{8}\times \Gamma_{8}:
  |\xi_1^*|\sim|\xi_2^*|\gtrsim N\gg |\xi_3^*|\}, \\
  B_2 =& \{(\xi_1,\cdots,\xi_{8},\tau_1,\cdots,\tau_{8})\in \Gamma_{8}\times \Gamma_{8}:
  |\xi_3^*|\gtrsim N\gg |\xi_4^*|\}, \\
  B_3 =& \{(\xi_1,\cdots,\xi_{8},\tau_1,\cdots,\tau_{8})\in \Gamma_{8}\times \Gamma_{8}:
  |\xi_4^*|\gtrsim N\}.
\endaligned\end{equation*}

\noindent\textbf{Estimate in $B_1$.} By  (\ref{EM8-2}) and
(\ref{EM8'-2}), we have
$$
|M_8+\widetilde{M}_8|\lesssim  {N_1^*}^{\frac{1}{2}} {N_3^*}^{\frac{1}{2}}.
$$
Therefore, similar to the estimate in $A_1$ in Proposition \ref{6-linear}, we have
\begin{equation*}\aligned
\mbox{ LHS of (\ref{P5.2-1})}\lesssim &\;
 N^{2s-2} \int_{B_1}
\frac{\widehat{f_1}(\xi_1,\tau_1)\cdots \widehat{f_8}(\xi_8,\tau_8)}
{|\xi_1^*|^{2s-\frac{1}{2}} \langle\xi_3^*\rangle^{\frac{1}{2}}\langle\xi_4^*\rangle\cdots\langle\xi_8^*\rangle}\\
\lesssim& \;
N^{-\frac{5}{2}+}\left\|I^{\frac{1}{2}-}_{\pm}(f_1^*,f_3^*)\right\|_{L^2_{xt}}
\left\|I^{\frac{1}{2}-}_{\pm}(f_2^*,f_4^*)\right\|_{L^2_{xt}}
\left\|J_x^{-1}f_5^*\right\|_{L^\infty_{xt}}\cdots
\left\|J_x^{-1}f_8^*\right\|_{L^\infty_{xt}}\\
\lesssim&\; N^{-\frac{5}{2}+}\> \|f_1\|_{Y_0^+}\|f_2\|_{Y_0^-}\cdots
\|f_7\|_{Y_0^+}\|f_8\|_{Y_0^-}.
\endaligned\end{equation*}

\noindent\textbf{Estimate in $B_2$.} By  (\ref{EM8-1}) and
(\ref{EM8'-1}), we have
\begin{equation}\label{EM88'}
|M_8+\widetilde{M}_8|\lesssim  N_1^*.
\end{equation}
Moreover, it satisfies that
$$
|\xi_1^*|-|\xi_3^*|\sim |\xi_1^*|\;\;\text{in}\;  B_2.
$$
Indeed, we have $|\xi_1^*|= |\xi_2^*|+ |\xi_3^*|+o(N_3^*)$ (see the proof of Lemma \ref{Sigma-Lem} for more details).
Therefore, similar to the estimate in $B_1$, we have
\begin{equation*}\aligned
\mbox{ LHS of (\ref{P5.2-1})}\lesssim &\;
 N^{3s-3} \int_{B_2}
\frac{\widehat{f_1}(\xi_1,\tau_1)\cdots \widehat{f_8}(\xi_8,\tau_8)}
{|\xi_1^*|^{2s-1} |\xi_3^*|^s\langle\xi_4^*\rangle\cdots\langle\xi_8^*\rangle}\\
\lesssim& \;
N^{-3+}\left\|I^{\frac{1}{2}-}_{\pm}(f_1^*,f_3^*)\right\|_{L^2_{xt}}
\left\|I^{\frac{1}{2}-}_{\pm}(f_2^*,f_4^*)\right\|_{L^2_{xt}}\\
&\;
\qquad\,\cdot
\left\|J_x^{-1}f_5^*\right\|_{L^\infty_{xt}}\cdots
\left\|J_x^{-1}f_8^*\right\|_{L^\infty_{xt}}\\
\lesssim&\; N^{-3+}\> \|f_1\|_{Y_0^+}\|f_2\|_{Y_0^-}\cdots
\|f_7\|_{Y_0^+}\|f_8\|_{Y_0^-}.
\endaligned\end{equation*}

\noindent\textbf{Estimate in $B_3$.} We only consider the worst
case: $|\xi_j|\gtrsim N$ for any $j=1,\cdots, 8$. By (\ref{EM88'})
and the similar estimates in $A_4$ in Proposition \ref{6-linear}, we
have
\begin{equation*}\aligned
\mbox{ LHS of (\ref{P5.2-1})}\lesssim&\; N^{8s-8}\displaystyle\int_{B_3}
\frac{\widehat{f_1}(\xi_1,\tau_1)\cdots \widehat{f_8}(\xi_8,\tau_8)}
{|\xi_1^*|^{2s-1}|\xi_3^*|^{s} \cdots|\xi_8^*|^{s}}\\
\lesssim&\;
N^{-6+}\left\|f_1^*\right\|_{L^{6-}_{xt}}\cdots
\left\|f_5^*\right\|_{L^{6-}_{xt}}
\left\|J_x^{0-}f_6^*\right\|_{L^{6+}_{xt}}
\left\|J_x^{-\frac{1}{2}-}f_7^*\right\|_{L^\infty_{xt}}
\left\|J_x^{-\frac{1}{2}-}f_8^*\right\|_{L^\infty_{xt}}\\
\lesssim&\; N^{-6+}\> \|f_1\|_{Y_0^+}\|f_2\|_{Y_0^-}\cdots
\|f_7\|_{Y_0^+}\|f_8\|_{Y_0^-},
\endaligned\end{equation*}
This gives
the proof of the proposition.
\end{proof}

\begin{prop} \label{10-linear}
For any $s\geq\frac{1}{2}$, we have
\begin{equation}
 \left|\displaystyle\int_0^\delta\Lambda_{10}(M_{10};w(t))\,dt\right|
 \lesssim
 N^{-3+}\>
 \|Iw\|_{Y_1}^{10}.
 \label{Lambda10}
\end{equation}
\end{prop}
\begin{proof} When $|\xi_1|,\cdots, |\xi_{10}|\ll N$, we have  $M_{10}=0$.
Therefore, we may assume that $|\xi_1^*|\sim |\xi_2^*|\gtrsim N$. Further, by symmetry, we may assume $|\xi_1|\geq \cdots
\geq|\xi_{10}|$ again. Similar to (\ref{P5.1-1}), it suffices to show
\begin{align}\label{P5.3-1}
&\displaystyle\int_{\Gamma_{10}\times \Gamma_{10}}\!\!
\frac{M_{10}(\xi_1,\cdots, \xi_{10})\widehat{f_1}(\xi_1,\tau_1)
\cdots \widehat{f_{10}}(\xi_{10},\tau_{10})} {\langle\xi_1\rangle
m(\xi_1)\cdots \langle\xi_{10}\rangle m(\xi_{10})}\\
&\qquad\qquad\qquad \lesssim
N^{-3+}\>\|f_1\|_{Y_0^+}\|f_2\|_{Y_0^-}\cdots
\|f_9\|_{Y_0^+}\|f_{10}\|_{Y_0^-},\nonumber
\end{align}
where $\Gamma_{10}\times \Gamma_{10}
=\{(\xi_1,\cdots,\xi_{10},\tau_1,\cdots,\tau_{10}):\xi_1+\cdots+\xi_{10}=0,
\tau_1+\cdots+\tau_{10}=0\}$. Now we divide it into two regions:
\begin{equation*}\aligned
  D_1 =& \{(\xi_1,\cdots,\xi_{10},\tau_1,\cdots,\tau_{10})\in \Gamma_{10}\times \Gamma_{10}:
  |\xi_2|\gtrsim N\gg |\xi_3|.\}, \\
  D_2 =& \{(\xi_1,\cdots,\xi_{10},\tau_1,\cdots,\tau_{10})\in \Gamma_{10}\times \Gamma_{10}:
  |\xi_3|\gtrsim N\}.
\endaligned\end{equation*}

\noindent\textbf{Estimate in $D_1$.} By Lemma \ref{Sigma-Lem}, we
have $|\sigma_6|\lesssim 1$ and thus
\begin{equation}\label{EM10}
|M_{10}|\lesssim 1.
\end{equation}
Similar to the estimates in $A_1$ in Proposition \ref{6-linear}, we
have
\begin{equation*}\aligned
\mbox{ LHS of (\ref{P5.3-1})}\lesssim &\;
 N^{2s-2} \int_{D_1}
\frac{\widehat{f_1}(\xi_1,\tau_1) \cdots
\widehat{f_{10}}(\xi_{10},\tau_{10})} {|\xi_1|^{s} |\xi_2|^{s}
\langle\xi_3\rangle \cdots \langle\xi_{10}\rangle}\\
\lesssim& \;
N^{-3+}\left\|I^{\frac{1}{2}-}_{-}(f_1,f_3)\right\|_{L^2_{xt}}
\left\|I^{\frac{1}{2}-}_{-}(f_2,f_4)\right\|_{L^2_{xt}}\\
&\;
\qquad\,\cdot
\left\|J_x^{-1}f_5\right\|_{L^\infty_{xt}}\cdots
\left\|J_x^{-1}f_{10}\right\|_{L^\infty_{xt}}\\
\lesssim&\; N^{-3+}\> \|f_1\|_{Y_0^+}\|f_2\|_{Y_0^-}\cdots
\|f_9\|_{Y_0^+}\|f_{10}\|_{Y_0^-}.
\endaligned\end{equation*}

\noindent\textbf{Estimate in $D_2$.} We only consider the worst
case: $|\xi_j|\gtrsim N$ for any $j=1\cdots 10$. Thus by
(\ref{EM10}), and the similar estimates in $B_3$ in Proposition
\ref{8-linear}, we have
\begin{equation*}\aligned
\mbox{ LHS of (\ref{P5.3-1})}\lesssim&\; N^{10s-10}\displaystyle\int_{D_2}
\frac{\widehat{f_1}(\xi_1,\tau_1)\cdots \widehat{f_{10}}(\xi_{10},\tau_{10})}
{|\xi_1|^{s}|\xi_2|^{s}|\xi_3|^{s} |\xi_4|^{s}\cdots|\xi_{10}|^{s}}\\
\lesssim&\;
N^{-8+}\left\|f_1^*\right\|_{L^{6-}_{xt}}\cdots
\left\|f_5^*\right\|_{L^{6-}_{xt}}
\left\|J_x^{0-}f_6^*\right\|_{L^{6+}_{xt}}\\
&\;
\cdot\left\|J_x^{-\frac{1}{2}-}f_7^*\right\|_{L^\infty_{xt}}\cdots
\left\|J_x^{-\frac{1}{2}-}f_{10}^*\right\|_{L^\infty_{xt}}\\
\lesssim&\; N^{-8+}\> \|f_1\|_{Y_0^+}\|f_2\|_{Y_0^-}\cdots
\|f_9\|_{Y_0^+}\|f_{10}\|_{Y_0^-}.
\endaligned\end{equation*}
This gives
the proof of the proposition.
\end{proof}

\section{A comparison between $E_I^1(w)$ and $E_I^3(w)$}
In this section, we show that the third generation modified energy
$E_I^3(w)$ is comparable to the first generation modified energy
$E_I^1(w)=E(Iw)$. In Section 5, we have shown that $E_I^3(w)$ is
almost conserved with a tiny increment. Then the result in this
section forecasts that $E_I^1(w)$ is also almost conserved with a
similar tiny increment (which will be realized in next section). Now
we state the result in this section.
\begin{lem} Let $s\geq \frac{1}{2}$,  then we have
\begin{equation}\label{E3-E1}
    \left|E_I^3(w(t))-E_I^1(w(t))\right|\lesssim N^{0-}\left(\|Iw(t)\|_{H^1}^4+\|Iw(t)\|_{H^1}^6\right).
\end{equation}
\end{lem}
\begin{proof} By (\ref{E1}), (\ref{E2}) and  (\ref{E3}), we have
\begin{equation*}\aligned
E_I^3(w(t))-E_I^1(w(t))=&\;\frac{1}{2}\Lambda_4\left(M_4(\xi_1,\xi_2,\xi_3,\xi_4)
-\frac{1}{2}\xi_{13}m_1m_2m_3m_4;w(t)\right)\\
&\;+
\Lambda_6(\sigma_6;w(t)).
\endaligned\end{equation*}
Therefore, it suffices to prove
\begin{equation}\label{L7.1-1}
    \left|\Lambda_4\left(M_4(\xi_1,\xi_2,\xi_3,\xi_4)
-\frac{1}{2}\xi_{13}m_1m_2m_3m_4;w(t)\right)\right|\lesssim N^{0-}\|Iw(t)\|_{H^1}^4,
\end{equation}
and
\begin{equation}\label{L7.1-2}
    \left|\Lambda_6(\sigma_6;w(t))\right|\lesssim N^{0-}\|Iw(t)\|_{H^1}^6.
\end{equation}
For (\ref{L7.1-1}), we refer to (32) in \cite{CKSTT-02-DNLS}. Now we turn to
prove (\ref{L7.1-2}). By Plancherel's identity, it suffices to show
\begin{equation}\label{L7.1-3}
\displaystyle\int_{\Gamma_{6}}\!\!
\frac{\sigma_6(\xi_1,\cdots,
\xi_{6})\widehat{f_1}(\xi_1,t) \cdots
\widehat{f_{6}}(\xi_{6},t)} {\langle\xi_1\rangle
m(\xi_1)\cdots \langle\xi_{6}\rangle m(\xi_{6})} \lesssim
 N^{0-}\>\|f_1(t)\|_{L^2_{x}}\cdots
\|f_{6}(t)\|_{L^2_{x}}.
\end{equation}
We may assume that $|\xi_1|\geq |\xi_2|\geq \cdots |\xi_6|$ by
symmetry. Since $\sigma_6=0$ when $|\xi_j|\ll N$ for any
$j=1,\cdots,6$, we may assume that $|\xi_1|\sim |\xi_2|\gtrsim N$.
By Lemma \ref{Sigma-Lem}, we have $|\sigma_6|\lesssim 1$. Note that
$$
\langle\xi\rangle m(\xi) \gtrsim \langle\xi\rangle^s, \mbox{ for
any } \xi\in \R,
$$
we have by Sobolev's inequality,
\begin{equation*}\aligned
\mbox{ LHS of (\ref{L7.1-3})}\lesssim&\; N^{-2+}\displaystyle\int_{\Gamma_{6}}
\frac{\widehat{f_1}(\xi_1,\tau_1)\cdots \widehat{f_{10}}(\xi_{10},\tau_{10})}
{\langle\xi_3\rangle^{s+}\cdots\langle\xi_6\rangle^{s+}}\\
\lesssim&\; N^{-2+}\left\|f_1(t)\right\|_{L^{2}_{x}}
\left\|f_2(t)\right\|_{L^{2}_{x}}
\left\|J_x^{-\frac{1}{2}-}f_3(t)\right\|_{L^\infty_{x}} \cdots
\left\|J_x^{-\frac{1}{2}-}f_{10}(t)\right\|_{L^\infty_{x}}\\
\lesssim&\; N^{-2+}\> \|f_1(t)\|_{L^{2}_{x}}\cdots
\|f_{10}(t)\|_{L^{2}_{x}}.
\endaligned\end{equation*}
This gives
the proof of the lemma.
\end{proof}

\section{The Proof of Theorem 1.1}

\subsection{A Variant Local Well-posedness} In this subsection, we
will establish a variant local well-posedness as follows.
\begin{prop} \label{prop:modified-local}
Let $s\geq \frac{1}{2}$, then Cauchy problem (\ref{eq:DNLS2}) is
locally well-posed for the initial data $w_0$ satisfying $Iw_0\in
H^1(\R)$. Moreover, the solution exists on the interval $[0,\delta]$
with the lifetime
\begin{equation}
 \delta\sim \|I_{N,s}w_0\|^{-\mu}_{H^1}\label{delta}
\end{equation}
for some $\mu>0$. Furthermore, the
solution satisfies the estimate
\begin{equation}
    \|I_{N,s}w\|_{Y_1}
    \lesssim
    \|I_{N,s}w_0\|_{H^1}.\label{LSE}
\end{equation}
\end{prop}
\begin{proof}
By the standard iteration argument (see cf. \cite{Ta-99-DNLS-LWP}),
it suffices to prove the multilinear estimates,
\begin{equation}\label{MLE1}
\left\|I(w_1\partial_x\overline{w_2}w_3)\right\|_{Z_1} \lesssim
\|Iw_1\|_{Y_1}\|Iw_2\|_{Y_1}\|Iw_3\|_{Y_1},
\end{equation}
and
\begin{equation}\label{MLE2}
\left\|I(w_1\overline{w_2}w_3\overline{w_4}w_5)\right\|_{Z_1}
\lesssim \|Iw_1\|_{Y_1}\cdots\|Iw_5\|_{Y_1}.
\end{equation}

By Lemma $12.1$ in
\cite{I-team:2004:multilinear-estimates-for-periodic-KdV-and-applications},
it suffices to prove the multilinear estimates,
\begin{equation}\label{MLE1}
\left\|w_1\partial_x\overline{w_2}w_3\right\|_{Z_s} \lesssim
\|w_1\|_{Y_s}\|w_2\|_{Y_s}\|w_3\|_{Y_s},
\end{equation}
and
\begin{equation}\label{MLE2}
\left\|w_1\overline{w_2}w_3\overline{w_4}w_5\right\|_{Z_s} \lesssim
\|w_1\|_{Y_s}\cdots\|w_5\|_{Y_s}.
\end{equation}
These were proved in \cite{Ta-99-DNLS-LWP}.
\end{proof}

\subsection{Rescaling}
We rescale the solution of (\ref{eq:DNLS2}) by writing
$$
w_\mu(x,t)=\mu^{-\frac{1}{2}}w(x/\mu,t/\mu^2);\quad
w_{0,\mu}(x)=\mu^{-\frac{1}{2}}w_0(x/\mu).
$$
Then $w_\mu(x,t)$ is still the solution of (\ref{eq:DNLS2}) with the
initial data $w(x,0)= w_{0,\mu}(x)$. Meanwhile, $w(x,t)$ exists
on $[0,T]$ if and only if $w_\mu(x,t)$ exists on
$[0,\mu^2T]$.

By $m(\xi)\leq 1$ and (\ref{mass-energy law}), we know that
\begin{equation*}
\|Iw_\mu(t)\|_{L^2_x}\leq \|w_\mu(t)\|_{L^2_x}
=\|w_{0,\mu}\|_{L^2_x}=\|w_{0}\|_{L^2_x}<\sqrt{2\pi}.
\end{equation*}
This together with (\ref{GN}) yields
\begin{equation}\label{Iw-energy}
\|\partial_xIu_\mu(t)\|_{L^2_x}^2 \sim E_I^1(w_\mu(t)),
\quad \|Iw_\mu(t)\|_{H^1_x}^2\lesssim E_I^1(w_\mu(t))+1.
\end{equation}
Moreover, by (\ref{E-I}), we get that
$$
\|\partial_xIw_{0,\mu}\|_{L^2} \lesssim N^{1-s}/\mu^s\cdot
\|w_0\|_{H^s}.
$$
Hence, if we choose $\mu\sim N^{\frac{1-s}{s}}$ suitably, we have $
\|Iw_{0,\mu}\|_{H^1}\leq 5. $ Thus we may take $\delta\sim
1$ by Proposition \ref{prop:modified-local}.

By standard limiting argument, the global well-posedness of $w$ in
$H^s(\R)$ follows if for any $T>0$, we have
\begin{equation*}
    \sup\limits_{0\leq t \leq T} \|w(t)\|_{H^s}\lesssim C(\|w_0\|_{H^s},T).
\end{equation*}
Further, in light of (\ref{E-I}) and (\ref{Iw-energy}), it suffices to show
\begin{equation}\label{B of Isol}
    \sup\limits_{0\leq t \leq \mu^2T} E_I^1(w_\mu(t))\lesssim C(T)
\end{equation}
for some  N. In the following subsection, we shall prove it by almost conservation law and iteration.

\subsection{Almost conservation law and iteration}
By (\ref{dE3}), we have
\begin{equation*}\aligned
E_I^3(w_{\mu}(t))=\;E_I^3(w_{0,\mu})+& \int_0^t\big(\Lambda_{6}(M_6\cdot \chi_{\Gamma_6\backslash \Omega};w(s))\; ds\\
+ &\int_0^t \Lambda_{8}(M_{8}+\widetilde{M}_{8};w(s))
+\Lambda_{10}(M_{10};w(s))\big)\,ds.
\endaligned\end{equation*}
By Proposition \ref{6-linear}--Proposition \ref{10-linear} and
(\ref{LSE}), we have for any $t\in [0,1]$,
\begin{equation*}\aligned
E_I^3(w_{\mu}(t))\leq &E_I^3(w_{0,\mu})+C_1N^{-\frac{5}{2}+}\left(
\|Iw_{\mu}\|_{Y_1}^{6}+\|Iw_{\mu}\|_{Y_1}^{8}+
\|Iw_{\mu}\|_{Y_1}^{10}\right)
\\
\leq &
E_I^3(w_{0,\mu})+C_2N^{-\frac{5}{2}+}.
\endaligned\end{equation*}
Thus,
\begin{equation*}\aligned
E_I^1(w_{\mu}(t))\leq &E_I^1(w_{0,\mu})+
\left(E_I^1(w_{\mu}(t))-E_I^3(w_{\mu}(t))\right)\\
&\;+\left(E_I^3(w_{0,\mu})-E_I^1(w_{0,\mu})\right)
+C_2N^{-\frac{5}{2}+}.
\endaligned\end{equation*}
Using (\ref{E3-E1}), choosing $N$ suitable large and applying the
bootstrap argument,  we obtain that for any $t\in [0,1]$,
$$
E_I^1(w_{\mu}(t))\leq 10.
$$
Repeating this process $M$ times, we obtain for any $t\in [0,M]$,
\begin{equation*}\aligned
E_I^1(w_{\mu}(t))\leq &E_I^1(w_{0,\mu})+
\left(E_I^1(w_{\mu}(t))-E_I^3(w_{\mu})\right)\\
&\;+\left(E_I^3(w_{0,\mu})-E_I^1(w_{0,\mu})\right)
+C_2MN^{-\frac{5}{2}+}.
\endaligned\end{equation*}
Therefore, by (\ref{E3-E1}) again, we have $E_I^1(w_{\mu}(t))\leq
10$ provided $M\lesssim N^{\frac{5}{2}-}$, which implies that the
solution $w_\mu$ exists on $[0,M\delta]\sim[0,N^{\frac{5}{2}-}]$.
Hence, $w$ exists on $[0,\mu^2T]$ with the relation
$$
N^{\frac{5}{2}-}\gtrsim \mu^2T\sim N^{\frac{2(1-s)}{s}}T.
$$
Thus we may take $T\sim N^{\frac{9s-4}{2s}-}$. When
$s\geq\frac{1}{2}$, we have $\frac{9s-4}{2s}>0$. This implies
(\ref{B of Isol}) by choosing sufficient large $N$, and thus
completes the proof of Theorem 1.1 .

\subsection*{Acknowledgements.} The authors were supported by the NSF of
China (No. 10725102, No. 10801015). The second author
was also partly supported by Beijing International Center for
Mathematical Research.

\end{document}